\documentclass[12pt]{article}
\usepackage{amsmath,amssymb,amsfonts,amsthm,mathrsfs}
\usepackage{kpfonts}
\usepackage[margin=1in]{geometry}
\usepackage{natbib}
\usepackage[breaklinks=true]{hyperref}
\usepackage{breakcites}
\usepackage{setspace}
\usepackage{color}
\usepackage{bbm}

\setcitestyle{authoryear,open={(},close={)}}
\nopagebreak
\hypersetup{
  colorlinks=true,
  linkcolor=blue,
  citecolor=red
}
{\theoremstyle{definition}

}
\newtheorem{theorem}{Theorem}
\newtheorem{lemma}{Lemma}
\newtheorem{claim}{Claim}

\newtheorem{proposition}{Proposition}
\newtheorem{corollary}{Corollary}



\newcommand{\argmax}{\operatornamewithlimits{argmax}}

\newcommand{\ind}[1]{{\mathbbm{1}_{\left\{{#1}\right\}}}}
\begin{document}

\title{Social Learning Equilibria\thanks{%
    We would like to thank the editor and four anonymous referees for helpful comments and suggestions. We thank seminar audiences in BI Norway, Caltech, Carnegie Mellon,
    Columbia, HEC Paris, INSEAD, Irvine, Madison, Moscow, Oxford, Paris School of Economics, Penn
    State and Stanford.  We would also like to thank Itai Arieli,
    Laura Doval, Federico Echenique, Luciano Pomatto, Pablo Schenone, Nicolas Vieille and
    Xavier Vives for illuminating discussions. Manuel Mueller-Frank
    gratefully acknowledges the financial support of the Spanish
    Ministry of Economy and Competitiveness
    (Ref. ECO2015-63711-P). Elchanan Mossel was supported by NSF
    DMS-1737944, ONR N00014-17-1-2598, and CCF 1665252. Omer Tamuz was
    supported by a grant from the Simons Foundation (\#419427).}}
\author{Elchanan Mossel\thanks{%
    MIT.}, Manuel Mueller-Frank\thanks{%
    IESE\ Business School.}, Allan Sly\thanks{%
    Princeton University.}\ and Omer Tamuz\thanks{%
    Caltech. }}
\maketitle

\begin{abstract}
We consider a large class of social learning models in which a group of agents face 
uncertainty regarding a state of the world, share
the same utility function, observe private signals,   and interact in a general dynamic  setting. We
introduce Social Learning Equilibria, a static  equilibrium concept that
abstracts away from the details of the  given extensive form, but
nevertheless captures the corresponding  asymptotic equilibrium behavior. We
establish general conditions for  agreement, herding, and information aggregation in
equilibrium,  highlighting a connection between agreement and information 
aggregation.

\noindent{\normalsize{{\bf Keywords}: Consensus; Information Aggregation; Herding.} }%
\end{abstract}
\section{Introduction}

Social learning refers to the inference individuals draw from observing the behavior of others, an inference which in turn impacts their own behavior. Social learning has served as an explanation for economic phenomena such as
herding\footnote{%
See \cite{banerjee, bhw, ss}.}, bubbles and crashes in financial markets%
\footnote{%
E.g., \cite{ss2, welch, ck}.} and technology adoption.\footnote{%
E.g., \cite{wb,duan}.} 

Many theoretical models of rational social learning are based on a
given dynamic social learning game, specified by a \emph{social learning
setting} and an \emph{extensive form.} The setting consists of the players,
their actions and common utility function, the state and signal spaces, and
a commonly known joint probability distribution over the state and signals. The extensive form
specifies the decision times of players and what each player observes at every given decision instant. These include sequential models, models of learning on social networks%
\footnote{E.g., \cite{parikh, gk03, rsv, adlo, mst, mmf, ls}. Another early contribution, in a boundedly rational framework, is due to \cite{bala1998learning}.} and more.

This approach has two inherent weaknesses. First, the analysis of asymptotic
equilibrium behavior in dynamic games is not straightforward\footnote{\cite{gk03}: ``The  computational difficulty of solving the model is massive even
in the  case of three persons [...] This is an important subject for future 
research.''}, resulting in a limited range of tractable models and a focus
on extremely stylized settings; in particular, the literature has largely avoided studying models of repeated actions by rational non-myopic players.\footnote{Exceptions include \cite{rsv} and \cite{mst}.} Second, when trying to understand or predict
behavior in ``real world'' social learning settings, the modeler might not know
the exact nature of interaction among individuals, and the sensitivity of the conclusions to each aspect of the extensive form
is often unclear.

To address these issues we introduce a static equilibrium approach which we
call \emph{social learning equilibrium} (SLE). We abstract away from the
extensive form dynamics, and focus directly on the asymptotic steady state
to which the game converges. For a given social learning setting, an SLE
includes a description of the information available to each agent, and an
action chosen by each player, as a function of what she knows. An SLE does
not include any details of the extensive form. However, the information
available to the agents was presumably learned through participation in some extensive
social learning game. The equilibrium condition simply requires each agent's
action to be optimal, given her information.

We consider large groups of agents and study three phenomena. {\em Agreement} occurs when all agents choose the same action.  {\em Herding}---a weaker form of agreement---occurs  when almost all agents choose the same action. Finally, {\em information aggregation} occurs when all agents choose the optimal action, given the realized state.

Our main results establish easily verifiable conditions on the setting and information structure which guarantee either agreement, herding or information aggregation. As our most important result, we point out a deep connection between herding and information aggregation: when private signals are unbounded,\footnote{As defined by \cite{ss}, private signals are unbounded if the support of the probability of either state conditional on one signal contains both zero and one.} herding can only occur when large amounts of information are exchanged, and in particular enough information must be exchanged for agents to learn the correct action. Perhaps surprisingly, this holds regardless of the extensive form, generalizing a result of \cite{ss} for the classical sequential model. The usefulness of all of our results lies in the fact that they are proved in the static SLE setting, making them robust to the details of the extensive form, and enabling the study of agreement, herding and information aggregation across a large spectrum of models---and in particular models of repeated actions---that are intractable to detailed dynamic analysis.



For most of our results we focus attention on the \emph{canonical setting}
of social learning with countably many agents, binary states and actions, a common prior and conditionally i.i.d.\ signals; this choice allows us to more easily explain our ideas and techniques. We also provide examples of how these can be applied beyond the canonical setting. 

\paragraph{A motivating example.} To more concretely describe our approach and results we study a particular social learning game, which a priori is not straightforward to analyze.  We explain how SLEs can be used to study this game, and what our results imply for it.

The setting of this game is the above mentioned {\em canonical setting}, with countably many agents, binary actions and states, and stage utilities which equal $1$ if the action matches the state, and $0$ otherwise.

The extensive form of this game is a variant of the classical sequential herding model, but with a repeated action twist. Agents are exogenously ordered. Agent $i$ chooses an action in each of the time periods $i, i+1, \ldots, i+100$, so that she lives for $101$ time periods, taking an action in each of these periods until her last action in period $i+100$, after which she leaves the game. Actions are public information, so that each agent observes all the actions taken by all the agents in all previous time periods, since period $1$.  Note that agents receive no additional information beyond their private signals and others' actions, and in particular observe their stage utilities only after taking their last action. Agents discount stage utilities by a common factor. 

To the best of our knowledge this particular game has not been previously studied. One obstacle to the analysis of this game is that in equilibrium, agents may not choose the myopically optimal action in each period, in an attempt to extract more information from future observations of their peers' actions.\footnote{Such equilibria have been constructed for other social learning games in \cite{mst}.} Indeed, it seems that one cannot hope to obtain a complete, detailed description of any equilibrium of this game, even if one considers myopic agents. Still, we now explain how our approach allows us to establish herding and information aggregation results for this model.

\paragraph{Social learning equilibria.} Instead of directly studying the dynamics of this game, we study the steady state reached at time infinity. This will be captured by an SLE, and thus will  consist of a (random) action for each player, as well as the information available to her. 

Given a Nash equilibrium of the extensive form game, the associated SLE consists of the last action of each player, and the information available to her at the end of the game. That is, for player $i$ it will include the action this player took at time $i+100$, and the information available to her when taking this action, which includes all the actions taken by others up to and including time $i+99$. Clearly, this action is a best response to this information, and so the equilibrium condition of an SLE is satisfied.

More generally, we show in Theorem~\ref{thm:game-equilibria} that also in games where agents take infinitely many actions---in fact, in a large class of such games---agents converge to an action (or actions) that is optimal given the information they have at the limit. Thus, every social learning game converges to an SLE.

\paragraph{Herding.} A natural question is whether a herd arises in this game, as in the classical sequential models. To study herding here, we will use only one additional feature of the extensive form. Namely, that each agent $i$ observes the actions of each agent $j < i$. Hence, in the corresponding SLE, the information available to agent $i$ includes the equilibrium actions of every agent $j < i$. We call such SLEs {\em weakly ordered}.

In Theorem~\ref{thm:herding} we show that every weakly ordered SLE satisfies {\em herding}: there is a (random) action that is taken by all agents, except finitely many. This immediately implies that in any equilibrium of the extensive form game described above, the last actions of all (but finitely many) of the agents are the same.

\paragraph{Information aggregation.} Our main result is Theorem~\ref{thm:learn-to-agree}. It establishes a fundamental connection between herding and  information aggregation.

Theorem~\ref{thm:learn-to-agree} states that when private signals are unbounded, information aggregation occurs in every  SLE that exhibits herding: the herding action is optimal, given the state. Thus, again, this applies to any equilibrium of the particular game described above. 
More generally, Theorem~\ref{thm:learn-to-agree} shows that---with unbounded signals and in
large groups---information aggregation is independent of the exact extensive
form, and even independent of the information agents have beyond their
private signals: it holds whenever herding occurs. Hence, this result shows that in order for agents to herd---that is, to almost all agree---they need to exchange a large amount of information, and in fact an amount large enough for them to herd on the optimal action.

\paragraph{Social networks and herding in probability.} One may be interested more generally in SLEs that are not weakly ordered. In particular, a large literature studies models of social learning on social networks, in which the observation structure is sparser. A variant of the game described above is one in which agent $i$ does not observe the actions of every agent $j < i$, but just those of some subset; such models (without the repeated action aspect) have been studied by \cite{adlo} and others. As an example of another application of our techniques, we consider the case that actions are {\em eventually public}: That is, the actions of agent $j$ are not observed by every $i > j$, but only by every $i$ {\em large enough}.  We call such SLEs {\em almost weakly ordered}.

As we explain below, such SLEs no longer necessarily satisfy herding. We show in Theorem~\ref{thm:co-finite} that they do, however, satisfy a weaker form of herding, which we call {\em herding in probability}: there is a (random) action that is taken with arbitrarily high probability by all the agents, except finitely many.\footnote{Formally, the agents actions converge {\em in probability} to some random herding action: for any $\varepsilon>0$ there are only finitely many agents $i$ whose action is unequal to the herding action with a probability greater than $\varepsilon$.} In our main result, Theorem~\ref{thm:learn-to-agree}, it in fact suffices to require herding in probability to induce information aggregation. Thus, also in this variant of our game, when signals are unbounded then all agents converge to the optimal action, with probability that tends to $1$ as $i$ tends to infinity.



\paragraph{Complete social learning equilibria.} Weakening the assumption of weakly ordered SLEs to almost weakly ordered SLEs resulted in the weaker result of herding in probability, rather than herding. It may also be interesting to {\em strengthen} the assumption, and see if we can derive a stronger conclusion.

A \emph{complete social learning equilibrium} (CSLE) is an SLE in which
each agent's information includes the other agents' equilibrium actions. Such SLEs arise, for example, in the limit of social learning games in which all agents act in all periods, and every agent observes the actions of all others \citep[e.g.,][]{sebenius1983don}. CSLEs are more generally related to models in which others' actions or beliefs are common knowledge, such as ``Agreeing to Disagree'' \citep{aumann}, the No Trade Theorem of \cite{milgrom-stokey}, and the common knowledge  equilibria of \cite{demarzo1998aggregation}.

Clearly, every CSLE is weakly ordered, and hence satisfies herding. Theorem \ref{thm:consensus} shows that, in a canonical setting, every CSLE in fact satisfies \emph{agreement}, i.e., all agents select the same action almost
surely, rather than just all but a finite number. Previous work
implies that agreement must hold unless agents are indifferent.\footnote{The initial insight is due to \cite{aumann}, with an important contribution by \cite{milgrom-stokey}. \cite{sg} studied the question for finite action settings. In the setting of social networks this was studied by \cite{mmf} and \cite{rsv}.} Our
contribution is to show that in large groups indifference is impossible, and
so agreement always holds.


These results imply that in any game that converges to a CSLE, all agents converge to the same action, and that when signals are unbounded they furthermore converge to the optimal action. This implication does not require any analysis of the extensive form, but merely that agents know which action their peers converge to. 
This includes countless possible models, 
some of which may potentially be intractable to detailed analysis.

\paragraph{Information aggregation in social networks.}
We apply our SLE approach to the repeated action model analyzed by \cite{mst}. In this social learning game all agents act in all periods. Agents are connected by a social network and only observe the actions of their neighbors.

We complement the results of \cite{mst}, showing in Proposition \ref{prop:mst-agreement} that asymptotic agreement holds in every Nash equilibrium, for any strongly connected network.
It then immediately follows that the agreement action satisfies information aggregation if private
signals are unbounded (Corollary \ref{cor:network-unbounded-learning}). In contrast, \cite{mst} study how network structure affects information aggregation, and provide sufficient conditions on the network that guarantee aggregation, regardless of whether signals are unbounded or not.

We also provide a new and simple proof of a result of \cite{mst}, showing that in symmetric networks the asymptotic agreement action satisfies information aggregation, for any informative signal structure. 

\paragraph{Concentration of Dependence.} A driving force behind our results is what we call the \emph{Concentration of Dependence Principle}. Informally, this principle refers to the fact that when an event $E$ is a function of i.i.d.\ random variables, then $E$ is approximately independent of almost all the random variables (Lemma~\ref{lem:principle}). Although we do not regard this principle as a novel contribution to probability theory, we believe that the value of applying it to economics might go beyond social learning applications.

In our social learning setting the Concentration of Dependence Principle implies that social learning outcomes in large groups depend on the state, and beyond that only on a small number of signals. We use this observation to prove almost all of our results: to preclude indifference (and hence disagreement) in complete SLEs, to show that herds arise in weakly ordered SLEs, and to show that herding implies information aggregation.

\subsection*{Extensions}

We consider several extensions of our model and results. Importantly, we study the case of bounded signals, where the support of the belief conditional on one signal contains neither zero nor one. Here we borrow the concept of \emph{information diffusion} introduced by \cite{ls} in the context of the sequential social learning model.\footnote{Assume that the support of private
beliefs is $[1-\beta ,\beta ]$. An action satisfies information diffusion if
it is optimal given the state with a probability of at least $1-\beta $.} We show that for bounded signals our theorems hold when one replaces information aggregation with the weaker notion of information diffusion. Other extensions include the cases of large finite groups of agents; games that include both rational and boundedly rational agents; settings with more than two states and actions; and agents with heterogeneous preferences and beliefs.\footnote{The last two extensions appear in the supplementary material.}

\subsection*{Related literature}

The social learning literature is too large to comprehensively cite here.\footnote{For a recent survey see \cite{gs}.} We
limit the discussion to those papers whose results are most closely related.

Our equilibrium approach is more in line with Aumann's approach \citeyearpar{aumann} of studying a static environment with common knowledge,
as compared to later social learning papers \citep[e.g.,][]{gp}, which
analyze the process by which common knowledge is reached. Similarly to
Aumann, we directly study the equilibrium, rather than specifying the exact
interaction structure and procedure by which the equilibrium is obtained.
Indeed, in many other fields of economics the tendency is to study static
equilibria directly rather than extensive forms. As we show in Theorem~\ref{thm:game-equilibria} there is no loss in restricting attention to SLEs when analyzing asymptotic equilibrium behavior of social learning games. 

Our results for the canonical setting provide new insights to two classes of social learning games that have been extensively analyzed in the literature; the repeated interaction model on social networks, and the canonical sequential social learning model. In particular, our Proposition~\ref{prop:mst-agreement} complements the agreement results for
settings of repeated interaction of \cite{gk03}, \cite{mmf}, and of \cite%
{rsv} which all show that agreement occurs but in case of indifference among
actions. 
Our Corollary ~\ref{cor:network-unbounded-learning} extends the results of \cite{mst} in showing that information aggregation holds for any strongly connected network  if signals are unbounded. 
In sequential settings, our work extends the classical herding and information aggregation results from sequential models  \citep{bhw, banerjee, ss} to a large class of social learning games.

\cite{celen} analyzed a variant of the canonical sequential social learning game where each agent only observes his immediate predecessor. \cite{adlo} and \cite{ls} extended this approach by considering the case where each agent observes a random subset of his predecessors.\footnote{\cite{adlo} assumes independent neighborhood draws across agents while \cite{ls} allow for correlation.} Our Theorem \ref{thm:herding} is motivated by their work, as it establishes a general sufficient condition on the observation structure that induces herding in probability.
Our Proposition \ref{prop:learn-to-agree-bounded} is closely related to \cite{ls}. For bounded
signals, they introduce the notion of information diffusion, which is a
weakening of information aggregation.
They provide two sufficient conditions on the random observation structure
such that information diffuses (respectively, fails to diffuse) in any
equilibrium. Applying our results to this concept, we shed additional
insight by connecting information diffusion to herding in probability.

The notion of an SLE is conceptually closely related to that of a correlated Bayes-Nash equilibrium. CSLEs are related to 
rational expectations equilibria. Theorem~\ref{thm:learn-to-agree} is similar to some results on optimality of rational expectations equilibria  \citep{dutta1997revelation,demarzo1998aggregation,dms}, which, unlike our results, already hold for a small number of players; a likewise similar approach is taken by \cite{ostrovsky2009information} in studying information aggregation in dynamic markets, and by \cite{babus2018trading} who study trading in over-the-counter markets.
We discuss these similarities and differences more thoroughly in \S \ref{sec:model}. \cite{minehart1999ex} introduce a concept of rational expectations equilibrium in a particular social learning setting. Despite some superficial similarities, their approach is essentially different from ours. For example, an equilibrium---as they define it---does not usually exist, and so they revert to an approximate
equilibrium notion, in which they prove their main results.

\bigskip

The rest of the paper is organized as follows. \S \ref{sec:model}
introduces the model and our equilibrium notion. \S \ref{sec:csle}
presents our results on agreement and information aggregation in CSLEs.
\S \ref{sec:herding} establishes our results on herding and
information diffusion in SLEs. \S \ref{sec:games} establishes the
formal relation between social learning equilibria and asymptotic
equilibrium behavior in social learning games. \S \ref{sec:networks} applies our results to models of social learning on networks. \S \ref%
{sec:extensions} presents some extensions.  \S \ref{sec:conclusion} is a conclusion.

\section{The Model}

\label{sec:model}

We consider a group of agents who must each choose an action under
uncertainty about a state of nature. Each agent's utility depends only on
her own action and the state, and agents are homogeneous in the sense of
sharing the same utility function and prior. Each agent observes a private signal, and
additionally some information about the others' signals. A social learning
equilibrium (SLE) includes a description of this additional information, as well as a  choice of action for each agent that maximizes her
expected utility, given the information available to her. We now define this formally.

\subsubsection*{Social learning settings}

A \emph{social learning setting} $(N,A,\Theta ,u,S,\mu )$ is defined by a
set of players $N\,$, a compact metrizable action space $A$, a
compact metrizable state space $\Theta $, a continuous utility function $%
u:A\times \Theta \rightarrow \mathbb{R}$, a measurable private signal space $%
S$, and finally a commonly known joint probability distribution $\mu $ over $%
\Theta \times S^{N}$.

We will denote by $\theta $ the random state of nature and by $\bar{s}%
=(s_{i})_{i\in N}$ the agents' private signals. When no ambiguity arises we
will denote probabilities and expectations with respect to $\mu $ by $%
\mathbb{P}[\cdot ]$ and $\mathbb{E}[\cdot ]$, respectively. For some
modeling applications it will furthermore be useful to add to this
probability space a non-atomic random variable $r$ that is independent of
the rest. This $r$ is an additional source of randomness that the players will use to implement mixed strategies.

\paragraph{Social Learning Equilibria (SLE).} Each agent $i$, in addition to her private signal $s_{i}$, learns $\ell _{i}$%
, which is some function of $\bar{s}$ (and possibly $r$). Agent $i$'s
(random) action is $a_{i}$. It takes values in $A$, and is some function of $%
\ell _{i}$ and $s_{i}$. Formally, $\ell _{i}$ and $a_{i}$ are random
variables that are, respectively, $\sigma (\bar{s},r)$- and $\sigma (\ell
_{i},s_{i})$-measurable.

Let $\bar{\ell}$ and $\bar{a}$ denote $(\ell_{i})_{i\in N}$ and $%
(a_{i})_{i\in N}$, respectively. In a given social learning setting, a \emph{%
\ social learning equilibrium} (or SLE) is a pair $(\bar{\ell},\bar{a})$
such that almost surely each agent's action $a_{i}$ is a best response,
given her information $\ell _{i}$ and $s_{i}$: 
\begin{equation}
a_{i}\in \operatornamewithlimits{argmax}_{a\in A}\mathbb{E}\left[ u(a,\theta
)\mid \ell _{i},s_{i}\right] \quad\quad\text{almost surely}.  \label{eq:sle}
\end{equation}

It is important to note that the information $\ell_i$ depends only on the private signals $\bar s$ and on the additional source of randomness $r$. In particular, $\ell_i$ does not depend explicitly on the actions $\bar a$, and so, if a player deviates, this has no effect on the others' information, and  thus there are no informational externalities. Therefore, since there are no payoff externalities, and since the action and state spaces are compact and utilities are continuous, an SLE exists for any setting. Moreover,  for any $\bar \ell$ there exists an $\bar a$---given by \eqref{eq:sle}--- such that $(\bar\ell,\bar a)$ is an SLE.

It is useful to think of $\ell_i$ as the information that player $i$ has
learned through participation in the equilibrium of some extensive form
game. Likewise, one should think of $a_i$ as the action that player $i$
converged to in the same game, so that the pair $(\bar{\ell},\bar{a})$
captures the asymptotic state of the game in question. We elaborate on this
in \S \ref{sec:games}, where we show that this asymptotic state
indeed satisfies the SLE condition~\eqref{eq:sle}.

So far we have put no restrictions on $\bar\ell$, and so, in this
generality, one would not expect to prove interesting results. In the
subsequent sections we will see how some relatively weak conditions on $\bar
\ell$ yield interesting properties of $\bar a$.

\paragraph{Complete Social Learning Equilibria (CSLE).} The first class of social learning equilibria which we consider are \emph{%
complete social learning equilibria} (or CSLE). In a CSLE each $\ell _{i}$ includes $\bar{a}$, in equilibrium. That is, the information available to each agent is sufficient to determine the equilibrium actions
of all other agents, and perhaps more information; formally, $a_i$ is $\sigma(\ell_j)$-measurable for every $i,j \in N$. We can therefore write $\ell_i = (\ell_i^0,\bar a)$ for some random variable $\ell_i^0$, and so, in equilibrium, it holds in a CSLE that
\begin{align}
\label{eq:csle}
a_{i}\in \operatornamewithlimits{argmax}_{a\in A}\mathbb{E}\left[ u(a,\theta
)\mid \ell_i^0,\bar{a},s_{i}\right] .  \end{align}

\paragraph{Weakly ordered and almost weakly ordered SLEs}

In a CSLE it holds for any two agents $i$ and $j$ that $a_i$ is $\sigma(\ell_j)$-measurable: in equilibrium, $j$ knows the action of $i$. This is a very strong assumption, which, as we show below, yields strong conclusions. We also study two weaker assumptions. 

We say that an SLE is \emph{weakly ordered} if the set of agents can be identified with the natural numbers  in such a way that if $i < j$ then
agent $j$ knows $i$'s equilibrium action: $a_{i}$ is $\sigma (\ell_{j})$-measurable.

An even weaker assumption is that of a {\em almost weakly ordered} SLE. An SLE is said to be almost weakly ordered if, for each agent $i$ there are only finitely many agents $j$ such that $a_i$ is not $\sigma(\ell_j)$-measurable. We show that these weaker assumptions imply  results which are weaker, but nevertheless have profound implications.

\paragraph{Discussion of the equilibrium concept.} 

SLEs can formally be thought of as a form of correlated (Bayes-Nash) equilibria, in which each agent receives a signal $\ell_i$ and best responds. Of course, since in the base game there are no externalities, this best response is independent of the actions of others. As is usual in correlated equilibria, we think of the signals $\bar \ell$ as being an endogenous part of the equilibrium \citep[see, e.g.,][\S 3.3]{or}. 
Mathematically equivalently, one could interpret the additional information $\bar\ell$ as being part of the environment, in which case the actions $\bar a$ are simply a Bayes-Nash equilibrium, as in \cite{bm}. We opt for the former interpretation of endogenous $\bar\ell$, as it better captures the position of an analyst who might not know $\bar\ell$ exactly. An additional motivation for this choice is our  application of SLEs to the study of the asymptotic state of an extensive form game. There, $\ell_i$ is the information learned by $i$ throughout the game, and in such games it is possible that agents' strategies affect what they learn.

CSLEs are related to rational expectations equilibria.\footnote{Similar ideas are also used in some versions of self-confirming equilibria \citep{rubinstein1994rationalizable, dekel1999payoff, dekel2004learning}. \cite{ostrovsky2009information} also uses a similar approach.} Our work is most closely related to that of \cite{demarzo1998aggregation}, and in particular the CSLE equilibrium condition \eqref{eq:csle} corresponds to condition (b) in their definition of a common knowledge equilibrium. Their results are also related to ours: in their setting, information is always aggregated optimally, whenever posterior estimates (which correspond to actions in CSLEs) are common knowledge. In our case this only holds for large groups and unbounded signals, as we show below.

\section{Agreement and Information Aggregation in Complete Social Learning
Equilibria}

\label{sec:csle}

In this section we study complete social learning equilibria (CSLEs). The CSLE assumption is strong, and thus yields strong, crisp conclusions; this makes CSLEs a good starting point for studying SLEs. In this section we also introduce the Concentration of Dependence Principle.

We focus on a class of social learning settings which appears frequently in the
literature: in \emph{canonical settings}, the set of players $N$ is countably infinite, $%
A=\Theta=\{0,1\}$, signals are informative and conditionally i.i.d., and $%
u(a,\theta)=1_{a=\theta}$, so that the utility is 1 when the action matches
the state, and 0 otherwise.

\paragraph{Agreement.} An SLE satisfies \emph{agreement} if almost surely $a_{i}=a_{j}$ for
all pairs of agents $i,j$. Our first result establishes agreement as a
property of any CSLE.

\begin{theorem}
\label{thm:consensus} In a canonical setting every CSLE satisfies agreement.
\end{theorem}

This result shows that Aumann's seminal agreement result carries over to
canonical social learning settings as a property of every CSLE. 
Previous results in the literature have
established that agreement is achieved, except in cases of indifference %
\citep{milgrom-stokey, rsv, mmf}. Our contribution is to show that, for the
case of CSLEs in canonical settings, indifference almost surely does not
occur and hence agreement holds. This essentially follows from what we call
the \emph{Concentration of Dependence Principle}, which we introduce now.
This principle underlies almost all of our results.

\paragraph{Concentration of Dependence.} Informally, concentration of dependence refers to the fact that when a
decision or event is a function of i.i.d.\ signals then it significantly
depends on only very few of them. The underlying mathematical fact is a well
known phenomenon known as \emph{mixing}, which we observe to have
interesting implications in our settings.

Formally, we say that a random variable $X$ is $\varepsilon $-independent of
an event $E$ if for every event $F$ that depends only on $X$ (i.e., if for
every $F\in \sigma (X)$) it holds that 
\begin{equation*}
\Big|\mathbb{P}[E\cap F]-\mathbb{P}[E]\cdot \mathbb{P}[F]\Big|<\varepsilon .
\end{equation*}%
Note that $X$ and $E$ are independent if and only if this holds for every $\varepsilon>0$.

\begin{lemma}[Concentration of Dependence Principle]
\label{lem:principle} Let $X_{1},X_{2},\ldots $ be independent random
variables, and let $E$ be any event defined on the same probability space.
Then except for at most $1/\varepsilon ^{2}$ many $i$'s, each $X_{i}$ is $%
\varepsilon $-independent of $E$.
\end{lemma}

For the convenience of the reader we provide a proof of this fact in \S %
\ref{app:concentration}.\footnote{%
Readers who are unfamiliar with this idea may wish to engage with some
examples. E.g., let $X_{1},\ldots ,X_{n}$ be i.i.d.\ fair coin tosses, and
consider two possible events. The first is the event that the majority of $%
X_{i}$'s equal $H$. It is easy to calculate and see that all the $X_{i}$'s
are very weakly correlated with $E$, and indeed intuitively this is clear. A
less obvious example is when $E$ is the event that an even number of $X_{i}$%
's are equal to $H$. Here, changing any $X_{i}$ (while keeping the rest
fixed) alters the indicator function of this event, and so it may seem that $%
E$ strongly depends on each $X_{i}$. However, $E$ is in fact independent of
each $X_{i}$.} In our canonical setting the private signals are i.i.d.,
conditional on the state. It thus follows from this principle that every
event that depends on the private signals is approximately conditionally
independent of almost all of them.

Let us now briefly outline the proof of Theorem~\ref{thm:consensus}. First
note that whenever both actions are taken in equilibrium, it must be that
all agents are indifferent between the actions. This follows from the same
intuition that underlies the no trade theorem of \cite{milgrom-stokey}, as
well as similar results in social learning 
\citep[e.g.,][]{sg, rsv,
  mmf}. This in turn implies that the disagreeing actions are equal to the
true state with probability $\frac{1}{2}$. Denote the disagreement event by $D$ and assume towards a contradiction that it has positive probability. For
agent $i$, let $b_{i}$ denote the optimal action of $i$ conditional only on
his private signal $s_{i}$. Note that since signals are informative, $b_{i}$
is equal to the state of the world with a probability strictly larger than $%
\frac{1}{2}$. Consider a deviation strategy of each agent $i$ such that she
selects the equilibrium action whenever $D$ does not occur and $b_{i}$
otherwise; this clearly cannot hurt the agent, since she is indifferent conditioned on $D$. The Concentration of Dependence Principle implies that $D$ is $\varepsilon$-independent of $s_{i}$ for all but finitely many agents $i$.
Hence for some agent $i$ (in fact, all but finitely many) the probability of 
$b_{i}$ being equal to the state conditional on $D$ is strictly larger than $%
\frac{1}{2}$, establishing a contradiction.

\paragraph{Information Aggregation.} We next turn to the learning properties of CSLEs. We have
shown above that the agents agree on the same (random)
action. Under which conditions is this agreement action optimal? That is, under which conditions is
an SLE information aggregating? Note that in a canonical setting an SLE is
\emph{information aggregating} if almost surely $a_{i}=\theta $ for all $i$.

To motivate the discussion we consider an example of a CSLE that is not information aggregating. In a canonical setting, assume that agents have a uniform prior, so that $\mathbb{P}[\theta=1]=\mathbb{P}[\theta=0]=1/2$, and that signals take values in $\{0,1\}$ with $\mathbb{P}[s_i=\theta \mid \theta] = 6/10$. 

For $i=1,\ldots,5$, let $\ell_i=(s_1,\ldots,s_5)$. That is, the first five players all know each others' private signals. Let $b$ be the optimal action given knowledge of the first five agents' signals:
$$
  b = \argmax_{a \in A}\mathbb{P}[a=\theta \mid s_1,s_2,s_3,s_4,s_5] = \begin{cases}
  1&\text{if }\sum_{i=1}^5 s_i \geq 3\\
  0&\text{if }\sum_{i=1}^5 s_i \leq 2.
  \end{cases}
$$

For $i>5$, let $\ell_i=b$. That is, for $i>5$, agent $i$'s additional information is the action $b$. For each $i$, let $a_i$ be some action satisfying the SLE condition
$$
  a_i \in \argmax_{a \in A}\mathbb{P}[a=\theta\mid \ell_i,s_i].
$$
Thus $(\bar\ell,\bar a)$ is an SLE.

As a simple calculation shows, in this case it holds that $a_i=b$ for all $i$: every agent is better off just following $b$ than doing anything else. Thus, this SLE is in fact a CSLE, as each agent knows $b$, and hence knows the actions of all other agents. Clearly, this CSLE is not information aggregating, since the probability that $b=\theta$ is not one.\footnote{In this same setting there exist other CSLEs that are information aggregating: for example, simply let all agents learn all the others' private signals. Then $a_i=\theta$ for all $i$, almost surely.}

The setting of this example is one with bounded signals; as defined by \cite{ss}, private signals are \emph{unbounded }if the support
of the {\em private belief} $p_{i}=\mathbb{P}[\theta =1\mid s_{i}]$ contains both $0$ and $1$. Similarly, private signals
are \emph{bounded} if the support of private beliefs contains neither $0$
nor $1$. The example above shows that with bounded signals a CSLE need not be information aggregating. The following result relates the
unbounded signal property to information aggregation in CSLEs.

\begin{proposition}
\label{prp:csle-aggregating} In a canonical setting with unbounded signals
every CSLE is information aggregating.
\end{proposition}

Combining Theorem \ref{thm:consensus} and Proposition \ref{prp:csle-aggregating}
we learn that in the case of CSLEs agreement implies information
aggregation, if signals are unbounded. This relation between agreement and information aggregation
holds much more generally, as we will establish in the next section. Indeed, Proposition~\ref{prp:csle-aggregating} is a corollary of a stronger result, Theorem~\ref{thm:learn-to-agree}. Nevertheless, we now sketch a proof of this Proposition, as it is simpler than that of Theorem~\ref{thm:learn-to-agree}.

The proof of Proposition \ref{prp:csle-aggregating} is also driven by the
Concentration of Dependence Principle. By Theorem~\ref{thm:consensus}, there
is some (random) agreement action $a_{0}$ that all players take. Consider
(towards a contradiction) the case in which the probability that $%
a_{0}=\theta $ is not 1, but some $q<1$. A player $i$ can consider the
deviation in which, instead of always choosing $a_{0}$, she chooses $a_{0}$
when her private signal is weak, but follows her private signal whenever her
private belief $p_{i}$ is strong. By strong we mean either greater than $q$
(in which case she would take action 1) or less than $1-q$ (in which case
she would take action 0). Because signals are unbounded, this occurs with
positive probability.

By the Concentration of Dependence Principle, $a_0$ is essentially a
function of some finite number of private signals, and almost all players $i$
have a private signal that is almost independent of the agreement action $%
a_0 $. Therefore this deviation is profitable for some (in fact, almost all)
players. Thus it is impossible that in equilibrium $q<1$, and so in
equilibrium $a_0=\theta$ almost surely.

\section{Herding and Information Aggregation}
\label{sec:herding}
Arguably, the most prominent result in the social learning literature is the
herding result established by \cite{bhw} in the canonical sequential social
learning model. They show that if agents make an irreversible binary
decision in strict sequential order, observing all the actions taken before
them, then eventually all agents take the same action. In other words,  \emph{herding} occurs: with probability one, all but a finite set of agents select the same action. The herding action is not necessarily optimal, even though the information contained in the pooled private signals suffices 
to choose the optimal action. \cite{ss} consider the canonical sequential social learning model and
show that when signals are unbounded herding still occurs, but the action chosen by the herd is optimal.

We analyze two different herding properties of SLEs.\footnote{The following definitions and equivalences apply for settings with a finite action set $A$.} We say that an SLE satisfies \emph{herding} if there is almost surely a cofinite set of agents who choose the same action. This (random) action is denoted as the \emph{herding action}. Equivalently, for any ordering of the agents, the sequence of random variables $(a_i)_i$ converges almost surely to the herding action $a^*$: 
\begin{align*}
\mathbb{P}[\lim_{i \to \infty}a_i=a^*]=1.
\end{align*}

We say that an SLE satisfies \emph{herding in probability} if there is a (random) herding action that the agents' actions converge to in probability. Formally, an SLE satisfies herding in probability if there is a random variable $a^*$ such that for any $\varepsilon>0$ there are only finitely many agents $i$ such that $\mathbb{P}[a_i \neq a^*] > \varepsilon$. Equivalently, for any ordering of the agents,
\begin{align*}
\lim_{i \to \infty} \mathbb{P}[a_i=a^*]=1.
\end{align*}

We would like to emphasize that despite the image that the term
``herd'' evokes, herding does not imply that the agents take a suboptimal action; indeed, the action
chosen by the herd can be correct with probability one \citep{ss}. Accordingly, we think of
herding as a weaker form of \emph{agreement}: an SLE satisfies agreement when all agents agree. A herding SLE is one in which
almost all the agents agree. Herding in probability holds when there is an action that only finitely many agents are significantly likely to disagree with. Note that agreement implies herding, which in turn implies herding in probability.

The subsequent analysis focuses on the relation between herding and information aggregation. Consider a SLE that satisfies either herding or herding in probability. We say that the herding action satisfies \emph{information aggregation} if it is almost surely optimal, conditioned on the state. In the canonical setting this means that the herding action is equal to $\theta$ with probability one.

\paragraph{Herding and information aggregation.}
Our first result of this section highlights a
deep connection between herding and information aggregation: when signals are unbounded, one cannot herd without
aggregating information. In other words, in order for agents to herd they must exchange a
large amount of information, and in particular an amount so large that they
learn the state in the process. 

\begin{theorem}
\label{thm:learn-to-agree} In a canonical setting, and when signals are
unbounded, in every SLE that satisfies herding in probability, the herding action satisfies information aggregation.
\end{theorem}

The proof of this theorem again relies on the Concentration of Dependence
Principle, and is similar to the proof of Proposition~\ref{prp:csle-aggregating}. The herding action is approximately independent of almost all the private
signals. Yet, it is taken by almost all the agents. Therefore, there will be
an agent who takes the herding action with very high probability, and whose
signal is almost completely independent from it. Hence such an agent would
prefer to follow her own private signal whenever doing so is more likely to
be correct than following the herd. But in equilibrium this agent does
follow the herd, and so it must be that her private signals never give an
indication that is stronger than the information contained in the herding
action. But this is impossible when signals are unbounded.

This result is related to similar results for rational expectations
equilibria 
\citep[e.g.,][]{dutta1997revelation,demarzo1998aggregation,dms,
  ostrovsky2009information}. There, however, agreement implies efficient
aggregation of information even for a small number of players and bounded
signals, whereas in our setting this holds less generally and crucially
depends on both the large size of the group and the unboundedness of the
signals.

\paragraph{A Sufficient Condition for Herding.} As we note above, herding is a weak form of agreement. Above we have shown that CSLEs satisfy agreement: if all agents observe each others' actions then they all agree. In this section we
relax the complete observation assumption of CSLEs to a weaker condition,
and show that it implies herding rather than agreement. 

Recall that an SLE is \emph{weakly ordered} if the set of agents can be identified with the natural numbers in such a way that if $i < j$ then
agent $j$ observes $i$'s action: $a_{i}$ is $\sigma (\ell _{j})$-measurable.

The classical sequential models of 
\cite{bhw} and \cite{ss} are a particular example in which the set of agents
is identified with the natural numbers and where each agent learns only the
actions of her predecessors. Weakly ordered SLEs are a larger class that
allows agents to furthermore have additional information, beyond the actions
of their predecessors. They can thus be used to model the asymptotic state
of games with richer extensive forms: perhaps the agents come in groups that
act together; perhaps they exchange information with the people standing
behind them or in front of them in line; or maybe they act more than once, as in the motivating example in the introduction of this paper.

We show that every weakly ordered SLE satisfies \emph{herding}

\begin{theorem}
\label{thm:herding} In a canonical setting every weakly ordered SLE
satisfies herding.
\end{theorem}

The proof of Theorem~\ref{thm:herding} uses similar ideas to that of Theorem~%
\ref{thm:consensus}, but involves a number of additional steps. Here,
one must first observe that if both actions are taken infinitely often then
agents must asymptotically be indifferent. If this occurs with positive
probability, then eventually agents will be able to guess (correctly with
high probability) that this will happen. Since---again
asymptotically---almost all agents have signals that are independent of this
event, they would choose to ignore it and follow their own private signals.
But then they would not be indifferent, and thus this cannot happen with
positive probability.


\paragraph{A Sufficient Condition for Herding in Probability.}
Consider an SLE where the $\ell _{i}$'s feature a layered observation structure: Assign each agent to one layer $(L_1,L_2,\ldots)$, with each layer having a finite number of agents. Assume that each agent observes, in addition to her private signal, the actions of all the agents in all the previous layers. Assume also private signals are unbounded \citep[see, e.g.,][]{rosenberg2017efficiency}.

Fixing the sizes of layers $L_1,\ldots,L_{n-1}$, there is, conditioned on the state, some non-zero probability that any agent in layer $L_n$ will choose action $0$, and some (other) non-zero probability that she will choose $1$. Since the actions of the players in a given layer are i.i.d.\ conditioned on the state and the previous players' actions, it follows from the Law of Large Numbers that if we choose each layer to be large enough, then with large probability there will, in every layer, be agents who choose both actions. Thus this SLE does not exhibit herding. However, we show that it does exhibit {\em herding in probability}. More generally, we show that herding in probability is obtained whenever each agent's action is ``eventually public''; that is, when each agent is observed by all except a finite group.

This is an example of an almost weakly ordered SLE. Recall that an SLE $(\bar \ell,\bar a)$ is almost weakly ordered if for each agent $i$ 
there are only finitely many other agents $j$ such that $a_i$ is not $\sigma(\ell_j)$-measurable. That is, the set of agents who observe agent $i$'s action is cofinite.
A natural example for an almost weakly ordered observation structure is the canonical sequential social learning model with the additional assumption that agents might not observe the actions of others who decided within a certain time interval before. However, if a given predecessor acted sufficiently earlier then her action is observed. One particular instance of this class is the above mentioned example of agents who are arranged in layers, and the agents of each layer observe the actions of all the agents in the previous layers. More generally, one could have a complicated social network structure in the spirit of \cite{adlo}.

\begin{theorem}
\label{thm:co-finite}
In a canonical setting, every almost weakly ordered SLE satisfies herding in probability.
\end{theorem}
It follows from this theorem and from Theorem~\ref{thm:learn-to-agree} that when signals are unbounded then almost weakly ordered SLEs satisfy information aggregation. This constitutes a strengthening of the learning theorem of \cite{ss} to a much larger class of extensive forms, and in particular to models of partial observations structures in the spirit of \cite{adlo}.

\section{Social Learning Equilibria and Social Learning Games}

\label{sec:games}

In this section we consider a large class of \emph{social learning games}. A
social learning game is a dynamic game with incomplete information in which
agents choose actions and observe information about other agents' actions
and signals. Its definition includes a \emph{social learning setting}---as
in the definition of an SLE---and a description of the extensive form. This
class comprises many models studied in the literature, including sequential
learning models, models of repeated interaction on social networks, and the game described in the introduction of this paper.

The main result of this section relates social learning games to SLEs. We
show that the asymptotic equilibrium behavior of agents in any social
learning game is captured by an SLE: for any distribution over asymptotic
equilibrium action profiles of a social learning game there exists an SLE
with a matching distribution over action profiles.

This correspondence provides motivation for studying SLEs, and also allows
to understand the long-run behavior of agents in many dynamic settings, by applying our results to the corresponding SLEs.

\subsubsection*{Social learning games}

A \emph{social learning game} includes a social learning setting $%
(N,A,\Theta ,u,S,\mu )$, together with a description of the extensive form
by which agents interact and learn. The extensive form consists of the tuple 
$(T,\mathfrak{k},\delta )$. For each agent $i$ the set $T_{i}\subseteq
\{1,2,\ldots \}$ denotes the set of \emph{action times} of agent $i$, i.e., the set
of time periods in which agent $i$ exogenously ``wakes
up'', receives information, and takes an action. The set $%
T=\left( T_{i}\right) _{i\in N}$ denotes the tuple of action times. For each
agent $i$ and time $t\in T_{i}$, let $\mathfrak{k}_{i,t}$ be the information
learned by agent $i$ at time $t$, and let $a_{i,t}$ be the action taken by
agent $i$ at time $t$. We denote by 
\begin{equation*}
\mathfrak{k}_{i}^{t}=\{\mathfrak{k}_{i,\tau }\,:\,\tau \leq t,\tau \in
T_{i}\}
\end{equation*}%
the information observed by agent $i$ by time $t$, and by 
\begin{equation*}
\mathfrak{k}_{i}=\{\mathfrak{k}_{i,t}\,:\,t\in T_{i}\}
\end{equation*}%
all the information observed by her, excluding her signal. We denote by 
\begin{equation*}
h_{i}^{t}=\{a_{i,\tau }\,:\,\tau <t,\tau \in T_{i}\}
\end{equation*}%
the actions taken by agent $i$ \emph{before} time $t$, and by 
\begin{equation*}
h^{t}=(h_{i}^{t})_{i\in N}
\end{equation*}%
all the actions taken by all the agents before time $t$.

The information $\mathfrak{k}_{i,t}$ is some function of the agents' actions
before time $t$, the private signals themselves, as well as the additional
independent random variable $r$, and takes values in some measurable space: 
\begin{equation*}
\mathfrak{k}_{i,t}=\mathfrak{k}_{i,t}(h^{t},\bar{s},r).
\end{equation*}%
The (possible) dependence on $r$ allows this framework to include mixed
strategies and random observation sets, such as observing a random subset of
the previously chosen actions.

The strategy of agent $i$ at time $t\in T_{i}$ is denoted by $\sigma _{i,t}$%
, takes values in $A$, and is some function of the information known to
agent $i$ at time $t$, which consists of $\mathfrak{k}_{i}^{t}$ and her
private signal $s_{i}$: 
\begin{equation*}
\sigma _{i,t}=\sigma _{i,t}(\mathfrak{k}_{i}^{t},s_{i}).
\end{equation*}%
The collection of maps $\sigma _{i}=(\sigma _{i,t})_{t\in T_{i}}$ is player $%
i$'s strategy, and the tuple of strategies across all agents, $(\sigma
_{i})_{i\in N}$, is the strategy profile. The history $\left( h^{t}\right)
_{t\in 
\mathbb{N}
}$ is generated according to $\left( \sigma _{i}\right) _{i\in N}$.

Finally, $\delta $ is the common discount factor, and agent $i$'s discounted
expected utility is 
\begin{equation*}
\sum_{t\in T_{i}}\delta ^{t}\cdot \mathbb{E}_{\sigma}[u(a_{i,t},\theta )].
\end{equation*}%
A strategy profile $\sigma $ is a \emph{Nash equilibrium} if for each agent $%
i$ her strategy $\sigma _{i}$ maximizes her discounted expected utility
given $\sigma _{-i}$, among all possible strategies for player $i$.

If agents are myopic, i.e. $\delta=0$, a strategy profile $\sigma$ is a Nash
equilibrium if for each agent $i$, given $\sigma_{-i}$, in each period $t$
her strategy $\sigma _{i,t}$ is such that her action $a_{i,t}$ maximizes her
expected utility in period $t \in T_i$ conditional on $\mathfrak{k}_i^t$ and 
$s_i$.

This definition of a social learning game is rather general and captures a
variety of different models. Most prominently it captures the sequential
social learning model of \cite{bhw}. To see this simply set $T_{i}=\{i\}$
for every agent $i$ and 
\begin{equation*}
\mathfrak{k}_{i,i}=\mathfrak{k}_{i}=\{a_{j,j}\,:\,j<i\}.
\end{equation*}

The sequential social learning models of \cite{adlo}, \cite{ls} and others
are likewise included in this framework. Here we have $T_{i}=\{i\}$ again,
but $\mathfrak{k}_{i,t}$ does not include all the actions of the
predecessors, but rather only those of a random subset of the predecessors.
The models of repeated interaction on social networks of \cite{gk03}, \cite%
{mossel2014asymptotic} and \cite{mst} can be captured by setting $T_{i}=%
\mathbb{N}$ for all agents $i$ and letting $\mathfrak{k}_{i,t}$ contain the
last period actions of all the neighbors of agent $i$. \cite{rsv} study a
more general model that is not subsumed by this framework, but still shares
many similarities. In fact, the proof of our result for this section,
Theorem~\ref{thm:game-equilibria}, exactly follows the proof of their
Proposition 2.1.

Finally, the models of repeated communication of beliefs in a social network
analyzed in \cite{gp} and \cite{parikh} can be captured by a
squared loss utility function and a discount factor equal to zero, hence
inducing myopic behavior.

For a given strategy $(\sigma _{i,t})_{t\in T_{i}}$, let $\bar{A}_{i}$
denote the (random) set of \emph{accumulation points} of agent $i$'s
realized actions; if $T_{i}$ is finite, then let $\bar{A}_{i}$ be the
singleton that contains only the last period action of agent $i$. If $T_{i}$
is infinite and $A$ finite, then $\bar{A}_{i}$ consists of the actions
chosen infinitely often.

Given these definitions, we are ready to establish the relation between Nash equilibria of social learning game and SLEs. As we mention above, this theorem is essentially due to \cite{rsv}.

\begin{theorem}
\label{thm:game-equilibria} Consider a social learning game, i.e., a social
learning setting and extensive form $(T,\mathfrak{k},\delta )$, and a
corresponding Nash equilibrium $\sigma $. Let $\bar{a}$ be any (random)
action such that $a_{i}\in \bar{A}_{i}$, and let $\ell _{i}=\mathfrak{k}_{i} 
$, where $\mathfrak{k}_{i}$ is generated according to $\sigma$. Then $(\bar{%
\ell},\bar{a})$ is an SLE.
\end{theorem}
This theorem  states that the asymptotic state of every Nash equilibrium\footnote{%
For this general class of social learning games there are a few possible definitions of a Perfect Bayesian equilibrium \citep[see e.g.,][]{watson2017general}, under all of which PBEs are Nash equilibria. Thus Theorem~\ref{thm:game-equilibria} applies to all PBEs.} is captured by an SLE.\footnote{The converse of Theorem~\ref{thm:game-equilibria} is also (trivially) true: given any SLE $(\bar \ell,\bar a)$, we can define the game in which $T_i=\{1\}$ for all $i$ and $\mathfrak{k}_{i,1}=\ell_i$ for each each $i$. Then there is an equilibrium in which every agent $i$ takes action $a_i$ at time $1$, and thus, under this equilibrium, this game (immediately) converges to the SLE $(\bar\ell,\bar a)$.} The information $\ell _{i}$ of agent $i$ is generated along the sequence of equilibrium actions in the social learning game, and thus $\ell _{i}$ depends on the interaction environment described by the extensive form. Importantly, $\ell_i$ depends also on the equilibrium strategies, since the information content of an agent's action depends on her strategy.

Theorem~\ref{thm:game-equilibria} essentially shows that the asymptotic equilibrium properties of \emph{any} social learning game can be analyzed via our SLE concept. Thus, while the definition of SLE is very permissive, its predictive power unfolds when applied to social learning games. 


In light of the definition of SLEs, Theorem~\ref{thm:game-equilibria} equivalently states that in any social learning game, every limit action of every agent $i$ is optimal conditional on her limit information $\mathfrak{k}_{i}$. This follows from Proposition 2.1 in \cite{rsv}. In the case that agent $i$ only acts finitely many times, it is immediate that her limit action---which, in this case, is by definition equivalent to her last action---is optimal conditioned on her information, and thus the SLE condition is satisfied.

The proof for agents that act infinitely often requires more work. The essential idea is that since beliefs converge, eventually each agent knows that her belief is unlikely to change substantially. Hence her incentive to deviate from the myopic expected utility maximizing action decreases.\footnote{This argument implies that although  we have assumed a common discount factor for all agents, Theorem ~\ref{thm:game-equilibria} still holds when discount factors differ across agents. }
We provide a version of the proof by \cite{rsv} adjusted to our language and notation in the appendix, establishing a link between their result and our concept of SLE. 

\subsubsection*{Learning and agreeing in social learning games}

In Theorem~\ref{thm:game-equilibria} we showed that the asymptotic outcomes
of social learning games correspond to SLEs. That is,
if for each agent $i$ we let $\ell_i$ denote the information $i$ has learned
by participating in the game, and if we let $a_i$ be the action that it
converged to (or some limit action in lieu of convergence), then $(\bar{\ell}%
,\bar{a})$ is an SLE. We can therefore apply our results on SLEs to social
learning games, with far-reaching implications. 

Consider any social learning game in a canonical setting in which each 
agent observes the limit action of all other agents; they could,
additionally, exchange information in other ways. Therefore, by Theorem~\ref{thm:game-equilibria}, the limit behavior of any Nash equilibrium of this game is captured by a CSLE. This implies, for example, that the repeated interaction model of \cite{gp} can be solved via CSLEs. It follows from Theorem~\ref{thm:consensus} that in the canonical setting all agents
must always converge to the same action. And if private signals are
unbounded, then by Proposition~\ref{prp:csle-aggregating} they must all converge
to the correct action.

A straightforward application of Theorem~\ref{thm:learn-to-agree} to social
learning games implies that in every social learning game in a canonical setting with unbounded signals, herding in probability
implies that
the herding action equals the realized state.

Theorem~\ref{thm:herding} implies that herding is indeed the outcome across
a large spectrum of social learning games: it suffices that if  $i < j$ then $%
j $ observes which actions $i$ converges to. This generalizes the results of 
\cite{bhw}, highlighting the deeper forces that drive them: herding (e.g.,
in the classical sequential model) is not a feature of the sequential timing
of actions, but rather of the observation structure of agents. In
particular, any social learning game with a weakly ordered observation
structure satisfies herding as an asymptotic equilibrium outcome; one such example is the game we introduce in the introduction of this paper. Relatedly, Theorem~\ref{thm:co-finite} shows that herding in probability holds in any social learning game where the limit action of every agent $i$ is observed by a cofinite set of agents.

\section{Social Learning in Networks}
\label{sec:networks}
In this section we use  SLEs and the Concentration of Dependence Principle to strengthen existing results on social learning in networks. More precisely, we apply our concepts to the social learning game analyzed by \cite{mst}. 
Denote the network neighbors of agent $i$ by $N_i$ and assume that the network is {\em strongly connected}: for each pair of agents $i$ and $j$ there is a finite tuple of agents $(k_1,k_2,\ldots,k_n)$ such that $k_1=i$, $k_n=j$, and $k_{m+1}$ is a neighbor of $k_m$ for $m=1,\ldots,n-1$. Each set of neighbors $N_i$ is assumed to be finite.

In the notation of \S \ref{sec:games}, their game can be described as follows: the setting is a canonical setting; the action times $T_i$ are the entire set $\mathbb{N}$ for every agent $i$; and the information observed by agent $i$ at time $t$ is the actions of her neighbors in the previous time period:
\begin{align*}
    \mathfrak{k}_{i,t} = (a_{j,t-1})_{j \in N_i}.
\end{align*}
Finally, utilities are discounted at a common rate. So agents observe their private signals in the beginning, at each period they take an action which yields stage utility 1 if it matches the (binary) state and 0 otherwise, and after taking this action they observe their neighbors' actions. Note that stage utilities are not observed, so that the initial private signals comprise all the available information.

As in \cite{mst} we assume that $\mathbb{P}[\theta=1]=1/2$, and that the private beliefs $\mathbb{P}[\theta=1\mid s_i]$ have a non-atomic distribution. Note that this in particular means that signals are informative.  As \cite{mst} show this implies that the sets $\bar A_i$ of accumulation points of the realized actions are almost surely the same for all agents. Namely, in every Nash equilibrium of this game there exists a (random) $\bar A$ such that $\mathbb{P}[\bar A_i=\bar A]=1$, for every agent $i$, so that when $\bar A = \{a\}$ then eventually each agent takes the action $a$, and when $\bar A=\{0,1\}$ then all agents choose both actions infinitely many times.

Using our Concentration of Dependence Principle in an argument identical to the one used in the proof of Theorem~\ref{thm:consensus}, it is possible to strengthen this result and show that $\bar A \neq \{0,1\}$. That is, the asymptotic indifference case cannot occur and hence all agents converge to the same (random) action $a^*$.
\begin{proposition}
\label{prop:mst-agreement}
Asymptotic agreement holds in every Nash equilibrium. That is, there is a (random) action $a^*$ such that
\begin{align*}
    \mathbb{P}[\bar A_i=\{a^*\}]=1
\end{align*}
for every agent $i$.
\end{proposition}

It follows from this proposition that for every Nash equilibrium of this game there is a unique SLE $(\bar\ell,\bar a)$ that is the limiting state of this game, in the sense of Theorem~\ref{thm:game-equilibria}: $\ell_i = \mathfrak{k}_i$ and $a_i \in \bar A_i = \{a^*\}$. Moreover, in this SLE we have agreement, since $a_i=a^*$. Hence  Theorem~\ref{thm:game-equilibria}, Theorem~\ref{thm:learn-to-agree} and Proposition~\ref{prop:mst-agreement} jointly imply the following corollary.
\begin{corollary}
\label{cor:network-unbounded-learning}
If signals are unbounded then in every Nash equilibrium the asymptotic agreement action satisfies information aggregation.
\end{corollary}

That is, the asymptotic agreement action $a^*$ is equal to $\theta$ almost surely. Different from Corollary~\ref{cor:network-unbounded-learning}, \cite{mst} focus on the properties of the network structure that assure information aggregation for all informative signal structures, rather than just unbounded ones. They sketch a proof of why this holds for symmetric equilibria on symmetric networks (e.g., the infinite two dimensional grid), and prove that, more generally, information is aggregated whenever no agent in the network is much more important than others.\footnote{More specifically, they show that information aggregation holds if there are numbers $d$ and $L$ such that no agent has observes more than $d$ others, and whenever and $i$ observes $j$, there is a path from $j$ back to $i$ of length at most $L$.} 

The formal proof of this result in \cite{mst} is rather involved and combinatorial. We provide here a short formal proof for the case of symmetric networks. We focus on the case that the network is simply the infinite chain, although the same proof applies more generally to symmetric networks.\footnote{By symmetric we mean that all agents play the same role in the geometry of the network. In the mathematics literature graphs with this property are called {\em vertex transitive}. Examples include  infinite grids and infinite regular trees. See, e.g., \cite{mst} for a formal definition.} Importantly, in this setting information aggregation is attained even when signals are bounded. We include the proof here, as it provides a short and illustrative example of an application of the Concentration of Dependence Principle.

For this proof we will need the additional assumption that each agent's strategy is some function of the entire vector of private signals, and does not depend on additional randomness. For pure equilibria this follows immediately from the definitions. In mixed equilibria this assumption is without loss of generality, since one can consider a model with equivalent outcomes in which we add to each signal an additional random component that is independent of the state, and which the agent can use to randomize, rather than using our global randomness $r$.

\begin{proposition}
\label{prop:mst-aggregation}
Assume the network is the infinite chain: agents are identified with the integers, and $i$ and $j$ are neighbors iff $|i-j|=1$. 
Then in any symmetric equilibrium the asymptotic agreement action $a^*$ satisfies information aggregation.
\end{proposition}
\begin{proof}
Condition on $\theta$. By the Concentration of Dependence Principle, all but finitely many private signals $s_i$ are $\varepsilon$-independent of $a^*$. By the symmetry assumption, if some $s_i$ is  $\varepsilon$-independent of $a^*$ then the same holds for all.  Hence all $s_i$'s are $\varepsilon$-independent of $a^*$. Since this holds for every $\varepsilon>0$, it follows that $a^*$ is independent of  each $s_i$.  By the same argument, $a^*$ is independent of the random variable $s_{[i,i+n]}=(s_i,s_{i+1},\ldots,s_{i+n})$ for any agent $i \in \mathbb{Z}$ and $n \geq 0$.

Note that $a^*$ is some function of the private signals, since it is the limit of any agent's actions, and each of these is some function of the private signals. Thus $a^*$ is simultaneously a function of the private signals and is conditionally independent of any finite set of them. Since signals are conditionally i.i.d., $a^*$ must depend on the tail of the sequence of private signals, and so must be constant, conditioned on the state. There are therefore four possibilities: either (1) $a^*=0$, or (2) $a^*=1$, or (3) $a^*=1-\theta$ or (4) $a^*=\theta$. Since in the first three cases the agents' utilities are at most $1/2$, these cannot be equilibria, as agents could profitably deviate by following their own private signals. Hence it must be that $a^*=\theta$.
\end{proof}
As \cite{mst} argue, an interpretation of this theorem is that egalitarianism leads to efficient aggregation of information. This is a message that has emerged, in various forms, in a number of other diverse settings \citep[e.g.,][]{kalai2004social, acemoglu2010spread, golub2010naive, dasaratha2017network}.\footnote{In other settings, \cite{dasaratha2018social} reach an opposing conclusion: diversity helps aggregation, while egalitarianism hurts it.} The intuition behind this result again stems immediately from the Concentration of Dependence Principle, which, combined with symmetry, starkly implies that actions are conditionally independent of private signals.

\section{Extensions}

\label{sec:extensions}

\subsection{Bounded signals}

Recall that Proposition~\ref{prp:csle-aggregating} shows that in a canonical
setting with unbounded signals every CSLE is information aggregating.

What can be said about information aggregation when signals are bounded?
Since independent of the signal structure there always exists an information
aggregating equilibrium,\footnote{%
To see this, consider the SLE where each agent's information includes the private signals of all others.} the question is what is the
worst possible equilibrium outcome in terms of learning. To answer this, we
borrow the notion of information diffusion introduced by \cite{ls} in
context of the sequential social learning model. Consider the support of the
private belief and let its convex hull be $\left[ \beta _{L},\beta _{H}%
\right] $. For simplicity assume that the support is symmetric, i.e., $\beta_L=\beta_H=\beta$.

A herding SLE in a canonical setting satisfies \emph{information diffusion} if the probability of the herding action $a^*$ being equal to the realized state $\theta$ is at least $1-\beta$. As Lobel and Sadler highlight, the notion of information diffusion is
particularly insightful if strong signals, i.e., those that induce a
posterior belief close to $\beta $ or $1-\beta $, are rare.



The next result is a generalization of Theorem~\ref%
{thm:learn-to-agree} to the bounded signal setting.

\begin{proposition}
\label{prop:learn-to-agree-bounded}In a canonical setting every SLE that
satisfies herding in probability also satisfies information diffusion.
\end{proposition}

\subsection{Large finite groups}

All of our theorems are proved in settings with infinitely many agents.
Analogous qualified statements for large finite groups follow from our
proofs.

For example, Theorem~\ref{thm:consensus} states that in a canonical
setting, every CSLE satisfies agreement, so that 
\begin{equation*}
\mathbb{P}[a_i=a_j \text{ for all } i,j]=1.
\end{equation*}


The following result is the analogous statement for agreement in CSLEs with finite groups of agents.

\begin{proposition}
\label{thm:consensus-finite} Consider a setting that is canonical, except that the group of agents is finite of size $n$. Then there is some constant $C>0$ that depends only on the distribution of private signals such that in every CSLE,
\begin{equation*}
\mathbb{P}[a_i=a_j \text{ for all } i,j]\geq 1-\frac{C}{\sqrt{n}}.
\end{equation*}
\end{proposition}
Proposition~\ref{thm:consensus-finite} shows that for a fixed private signal distribution, the probability of disagreement is at most of order $1/\sqrt{n}$, uniformly among all CSLEs with $n$ agents. The constant $C$, which is calculated explicitly in the proof, decreases as signals get more informative.

\subsection{Social Learning Games with Rational and Boundedly Rational Agents}
Earlier work in the literature on repeated interaction in social networks considered the case of social networks where rational and boundedly rational agents coexist. \cite{mmf2} shows that all agents asymptotically aggregate information if the set of actions is rich. \cite{clx} show that in the case of binary actions information aggregation can fail for certain network structures. 

It is easily verified that all asymptotic properties we establish via the SLE approach carry forward for the {\em subset of Bayesian agents in the network}, if the conditions we lay out in the theorems are satisfied for the Bayesian agents; the boundedly rational agents simply provide additional information to the rational ones, and this additional information can be embedded in $\ell_i$. Whether the results carry forward to the asymptotic behavior of the boundedly rational agents depends on their updating heuristic, and is beyond the scope of this paper. 


\section{Conclusion}

\label{sec:conclusion} We introduce social learning
equilibria as a useful tool to analyze social learning. 
We provide agreement, herding and information aggregation results for social
learning equilibria that unify and shed additional insight on the social
learning literature. In particular, we show that the relation between
unbounded signals and the optimality of the herding action established by 
\cite{ss} holds much more generally. In fact, in any canonical social learning
environment with unbounded signals the action selected in a Bayesian herd is
optimal.


There are several natural avenues for future research. The Concentration of Dependence Principle naturally lends itself to proving positive results: namely, that agreement, various forms of herding, and information aggregation occur when the agents observe enough of their peers' actions. These conditions on the observation structure are in general not necessary: indeed, in many cases it may be that a detailed analysis of the dynamics of an extensive form game can yield stronger conclusions than one can hope to deduce using the SLE approach.\footnote{One such example is the main result of \cite{mst}; we elaborate on this in \S\ref{sec:networks}.} We leave the pursuit of negative results and the associated necessary conditions to future work.





The second avenue of future research concerns the extension of our analysis to capture payoff externalities.  SLEs may prove to be  a useful tool to answer these questions, which we leave for future work.

\bigskip

\bibliographystyle{plainnat}
\bibliography{myrefs}

\appendix
\section{Concentration of Dependence}

\label{app:concentration} 
\begin{proof}[Proof of Lemma~\ref{lem:principle}]
  Choose any $\varepsilon > 0$ and let
  $X_{i_1},X_{i_2},\ldots,X_{i_k}$ be $k$ random variables that are
  not $\varepsilon$-independent of $E$.  Without loss of generality we
  may assume that $(i_1,i_2,\ldots,i_k) = (1,\ldots,k)$. Let
  $F_1,\ldots,F_k$ be events that witness the violation of
  $\varepsilon$-independence, so that for $i=1,\ldots,k$
  \begin{align}
    \label{eq:XAYB}
    \Big|\mathbb{P}[E \cap F_i]-\mathbb{P}[E]\cdot\mathbb{P}[F_i]\Big| \geq \varepsilon,
  \end{align}
  and each $F_i$ is in $\sigma(X_i)$.

  Let $Y$ be the indicator of the event $E$.  Then we can write
  \eqref{eq:XAYB} as
  \begin{align*}
    \Big|\mathrm{Cov}(F_i,Y)\Big| \geq \varepsilon.
  \end{align*}

  Let $\eta_i \in \{-1,+1\}$ equal the sign of
  $\mathrm{Cov}(F_i, Y)$. Then \eqref{eq:XAYB} is equivalent to
  \begin{align*}
    \mathrm{Cov}(\eta_iF_i,Y) \geq \varepsilon.
  \end{align*}
  Summing over $i$ we get
  \begin{align*}
    \sum_{i=1}^k\mathrm{Cov}(\eta_iF_i,Y) \geq k\varepsilon.
  \end{align*}
  By additivity of covariance, it follows that 
  \begin{align*}
    \mathrm{Cov}\left(\sum_{i=1}^k\eta_iF_i,Y\right) \geq k\varepsilon.
  \end{align*}
  By the Cauchy-Schwarz inequality it follows that
  \begin{align*}
    \sqrt{\mathrm{Var}\left(\sum_{i=1}^k\eta_iF_i\right)\cdot\mathrm{Var}(Y)} \geq k\varepsilon.
  \end{align*}
  Denote $I = \sum_{i=1}^k\eta_iF_i$, and note that
  $\mathrm{Var}(Y)\leq 1$, since $Y \in \{0,1\}$. So, squaring both
  sides yields $\mathrm{Var}(I) \geq k^2\varepsilon^2$. Note that $\mathrm{Var}(I)$ is at most
  $k$, since $I$ is the variance of $k$ independent random variables,
  each with variance at most 1. Hence we have that $k \geq k^2\varepsilon^2$, or $k \leq 1/\varepsilon^2$.
\end{proof}

\section{Proof of Theorem~\ref{thm:consensus} and Proposition~\ref{thm:consensus-finite}}

\label{sec:complete-proofs}  We start with the following lemma, which is essentially a
formulation of the No Trade Theorem of \cite{milgrom-stokey}. This lemma
states that when there is disagreement then players must be indifferent.

Let $(\bar \ell, \bar a)$ be a CSLE in a setting in which $\Theta=\{0,1\}$. Denote agent $i$'s equilibrium belief by $q_i = \mathbb{P}[\theta=1 \mid \ell_i, s_i]$,
and let the {\em disagreement event} $D$ be the event that $a_i\neq a_j$ for some $i,j \in N$.

\begin{lemma}
\label{lem:no-bet} Let $A=\Theta =\{0,1\}$ and $u(a,\theta )=1_{a=\theta }$
(as in a canonical setting, but with no restrictions on the signals). In any
CSLE, if the disagreement event $D$ has positive probability, then conditioned on $D$ it
almost surely holds that $q_{i}=1/2$ for all $i$.
\end{lemma}

\begin{proof}
  Consider an outside observer who observes all the agents' actions
  $\bar a$.  Her belief is $q_* = \mathbb{P}[\theta=1 \mid \bar a]$.
  Since $\bar a$ is $\sigma(\ell_i,s_i)$-measurable, it follows from
  the law of total expectations that for every $i$
  \begin{align}  \label{eq:q*}
    q_* = \mathbb{E}[q_i \mid \bar a].
  \end{align}
  Since $1$ is the action that is optimal for beliefs above $1/2$, we
  have that $a_i = 1$ implies that $q_i \geq 1/2$. Likewise, $a_i=0$
  implies $q_i \leq 1/2$. Hence the claim follows by~\eqref{eq:q*}.
\end{proof}

To prove Theorem~\ref{thm:consensus} we show that the probability of the disagreement event $D$ is zero. As we  show, this follows from the Concentration of Dependence Principle. 

The proof of this theorem will follow a strategy that we will use again for other results of this paper. We will consider, for each agent $i$, a deviation in which she plays a different action $b_i$ whenever she observes that the event $D$ occured, and otherwise plays $a_i$.

To define the action $b_i$ that she takes when she observes $D$, we consider (as an auxiliary construction) an additional  fictitious player $x$ who observes whether or not $D$ occurred (i.e., observes the random variable $\ind{D}$, the indicator of $D$), and additionally receives a signal $s_x$, which, conditioned on $\theta$, is  independent of the other agents' signals, and distributed identically to theirs. Denote by $b_x$  an action chosen by such a player, which is optimal conditioned on $D$:
  \begin{align}
  \label{eq:b1}
    b_x  =b(s_x) \in \operatornamewithlimits{argmax}_{a \in A}\mathbb{P}[a=\theta\mid D,s_x]. 
  \end{align}
The event $D$ depends on the (real) agents' actions, and so is conditionally independent of $s_x$. Note that we are not formally changing our model by adding agent $x$, but only using it as a way to define the function $b$ above.

Consider now a possible deviation by agent $i$ who chooses the action
\begin{align}
\label{eq:b2}
    b_i =b(s_i)
\end{align}
whenever she observes $D$. Here $b$ is the function defined in~\eqref{eq:b1}. Intuitively, agent $i$, by choosing $b_i$ when $D$ occurs, is choosing an action that would be optimal if $D$ were conditionally independent of her signal. In the following lemma we show that the probability that $b_i$ is the correct action tends to the probability that $b_x$ is the correct action. This follows from Concentration of Dependence. We state and prove this lemma in more generality than we need for this theorem, because we use it in other proofs.
\begin{lemma}
\label{lem:apply-cod}
  Let $(\bar\ell,\bar a)$ be an SLE in a canonical setting, and let $F$ be any event that is $\sigma(\bar a)$-measurable. Let
  \begin{align*}
    c_x  =c(s_x) \in \operatornamewithlimits{argmax}_{a \in A}\mathbb{P}[a=\theta\mid F,s_x],
  \end{align*}
  and let $c_i=c(s_i)$. Then
  \begin{align*}
    \lim_{i \to \infty}\mathbb{P}[c_i = \theta, F] = \mathbb{P}[c_x = \theta, F].
  \end{align*}
\end{lemma}
That is, for high $i$, under the event $F$, the probability that the action $c_i$ is correct tends to that of $c_x$, which---crucially---is conditionally independent of $F$.  

\begin{proof}[Proof of Lemma~\ref{lem:apply-cod}]
By the Concentration of Dependence Principle, for every $\varepsilon>0$ it holds for all $i$ large enough that both  $\theta=0$ and $\theta=1$
  \begin{align*}
      \big|\mathbb{P}[c_i = \theta, F \mid \theta] - \mathbb{P}
    [c_i = \theta\mid \theta]\cdot \mathbb{P}[F \mid \theta]\big| < \varepsilon,
  \end{align*}
  and so
  \begin{align*}
    \lim_{i \to \infty}\mathbb{P}[c_i = \theta, F \mid \theta] = \mathbb{P}%
    [c_i = \theta\mid \theta]\cdot \mathbb{P}[F \mid \theta],
  \end{align*}
  where the right-hand side does not depend on $i$, since the $c_i$'s are
  identically distributed. Since $c_x$ is also identically distributed it follows that
  \begin{align*}
    \lim_{i \to \infty}\mathbb{P}[c_i = \theta, F \mid \theta] = \mathbb{P}%
    [c_x = \theta\mid \theta]\cdot \mathbb{P}[F \mid \theta] = \mathbb{P}    [c_x = \theta,F\mid \theta],
  \end{align*}
  where the second equality holds since  $c_x$ and $F$ are conditionally independent. Multiplying both sides by $\mathbb{P}[\theta]$ and summing over $\theta=0$ and $\theta=1$ yields that also unconditionally
  \begin{align*}
    \lim_{i \to \infty}\mathbb{P}[c_i = \theta, F]  = \mathbb{P}    [c_x = \theta,F].
  \end{align*}
\end{proof}  
  
Given Lemma~\ref{lem:apply-cod} we are ready to prove our Theorem.

\begin{proof}[Proof of Theorem~\ref{thm:consensus} and Proposition~\ref{thm:consensus-finite}]
  To prove our claim we need to show that the probability of the disagreement event $D$ is zero.
  
  Assume towards a contradiction that $D$ has positive
  probability. We consider, for each player $i$, the  deviation
  of following her private signal whenever $D$ occurs, in which case she
  chooses action $b_i$ as in \eqref{eq:b2}. When the complement of $D$ occurs, she does not deviate, choosing $a_i$. 
  
  The profit player $i$
  stands to gain from this deviation is
  \begin{align*}
    P = \mathbb{P}[b_i=\theta,D]-\mathbb{P}[a_i=\theta,D].
  \end{align*}
  The second term is equal to
  \begin{align*}
      \mathbb{P}[a_i=\theta\mid D] \cdot \mathbb{P}[D] = \frac{1}{2}\mathbb{P}[D]
  \end{align*}
  since, by Lemma~\ref{lem:no-bet}, conditioned on $D$ player $i$ is indifferent, and so her
  expected utility from any action is $1/2$.
  
  Consider first the case that there are infinitely many agents (Theorem~\ref{thm:consensus}). By Lemma~\ref{lem:apply-cod}, the first term tends to 
  \begin{align*}
      \mathbb{P}[b_x=\theta,D] = \mathbb{P}[b_x=\theta\mid D] \cdot \mathbb{P}[D].
  \end{align*}
  Now, since signals are informative, and since conditioned on $D$ the probability of each state is $1/2$, $\mathbb{P}[b_x=\theta\mid D] > 1/2$. Hence  $P$ is positive for all $i$ large enough, and we have reached a contradiction with our equilibrium assumption. This completes the proof of Theorem~\ref{thm:consensus}.

  Consider now the case that there are $n < \infty$ agents (Proposition~\ref{thm:consensus-finite}).
  Denote $Q = \mathbb{P}[b_x = \theta]$.
  Since $b_x$ is conditionally independent of $D$, it follows that $Q = \mathbb{P}[b_x = \theta \mid D]$,
  and so $Q$ is simply the expected utility of agent $x$, conditioned on the disagreement event $D$, under which $x$ has posterior $1/2$. Therefore, since private signals are informative, $Q > 1/2$, and furthermore $Q$ depends on the private signal distributions, but not on the SLE: it is simply the expected utility of an agent whose prior makes her indifferent and who acts after observing a signal.
  
  Let $\varepsilon = 1/\sqrt{n/2-1}$. By the Concentration of Dependence Principle there are at least $n-1/\varepsilon^2 = n/2+1$ agents $i$ such that $b_i$ is $\varepsilon$-independent of $D$, conditioned on $\theta=1$, and likewise for $\theta=0$. Hence there is at least one agent $i$ such that $b_i$ is $1/\sqrt{n/2-1}$-independent of $D$ conditioned on both $\theta=1$ and $\theta=0$. Therefore, by the same calculation as in Lemma~\ref{lem:apply-cod}, it will hold for this agent that
  \begin{align*}
    \big|\mathbb{P}[b_i = \theta, D] - \mathbb{P}    [b_i = \theta]\cdot \mathbb{P}[D]\big|<\frac{1}{\sqrt{n/2-1}}.
  \end{align*}
  Since $\mathbb{P}[b_i=\theta]=\mathbb{P}[b_x=\theta]$, it follows that
  \begin{align*}
    \mathbb{P}[b_i = \theta, D] \geq Q \cdot
    \mathbb{P}[D] -\frac{1}{\sqrt{n/2-1}}.
  \end{align*}
  By the equilibrium assumption
    \begin{align*}
    \mathbb{P}[a_i = \theta, D] \geq \mathbb{P}[b_i = \theta, D],
  \end{align*}
  and so, since $\mathbb{P}[a_i=\theta,D]=\frac{1}{2}\mathbb{P}[D]$, we have that
  \begin{align*}
      \frac{1}{2}\mathbb{P}[D] \geq Q \cdot
    \mathbb{P}[D] -\frac{1}{\sqrt{n/2-1}}. 
  \end{align*}
  Rearranging yields
  \begin{align*}
      \frac{1}{Q-1/2}\frac{1}{\sqrt{n/2-1}} \geq \mathbb{P}[D].
  \end{align*}
  This completes the proof of Proposition~\ref{thm:consensus-finite}, for $C=\frac{1}{2(Q-1/2)}$.
\end{proof}

\section{Proof of Theorem~\ref{thm:learn-to-agree} and Proposition~\ref{prop:learn-to-agree-bounded}}

In this section we prove Proposition~\ref{prop:learn-to-agree-bounded}, which is a generalization of Theorem~\ref{thm:learn-to-agree}; Proposition~\ref{prop:learn-to-agree-bounded} applies more generally to signals that can be either bounded or unbounded. Since information diffusion coincides with information aggregation when signals are unbounded, Theorem~\ref{thm:learn-to-agree} follows immediately.

Let $(\bar\ell,\bar a)$ be an SLE that satisfies herding in probability. Let $a^*$ be the herding action, and let $p = \mathbb{P}[a^* = \theta]$
be the probability that the herding action is optimal. We need to show that $p=1$.

Let the convex hull of the support of private beliefs be $[\beta,1-\beta]$ (so that $\beta=0$ when signals are unbounded), and assume by contradiction that  $p \leq 1-\beta-2\varepsilon$ for some $\varepsilon > 0$. Hence, either for $\theta=0$ or for $\theta=1$ it holds that $\mathbb{P}[a^* = \theta \mid \theta] \leq 1-\beta-2\varepsilon$. Assume without loss of generality that this holds for $\theta=1$, so that
\begin{align*}
    p_1 = \mathbb{P}_1[a^* = \theta] \leq 1-\beta-2\varepsilon,
\end{align*}
where here, and in the remainder of this proof, we simplify notation by writing $\mathbb{P}_1[\cdot]$ to denote $\mathbb{P}[\cdot\mid \theta=1]$.

Similarly to the proof of Theorem~\ref{thm:consensus}, let 
\begin{align*}
b_i = b(s_i) \in \operatornamewithlimits{argmax}_{a \in A}\mathbb{P}%
[\theta=a \mid s_i]
\end{align*}
be an optimal action chosen given agent $i$'s private signal only.  Let $B_i$ be the
event that $\mathbb{P}[b_i=\theta\mid s_i] > 1-\beta-\varepsilon$. Since the 
$b_i$'s are identically distributed, all of the events $B_i$ have the same
probability. Furthermore, this probability is positive, by our assumption on
the support of the private signals.

Imagine that agent $i$ deviates and chooses $b_i$ whenever $B_i$ occurs, and
otherwise follows $a_i$. Then her expected gain in utility is 
\begin{align*}
\mathbb{P}[b_i=\theta, B_i] - \mathbb{P}[a_i = \theta,B_i].
\end{align*}
We prove the claim by showing that this gain is strictly positive, which contradicts the equilibrium assumption. In fact, we show that this already holds conditioned on $\theta=1$ (with the case $\theta=0$ following from the same argument). Recall that to simply notation we write $\mathbb{P}_1[\cdot]$ to denote $\mathbb{P}[\cdot\mid \theta=1]$.

To bound the first term, we note that, by the definition of $B_i$, 
\begin{align*}
\mathbb{P}_1[b_i=\theta, B_i] \geq (1-\beta-\varepsilon)\mathbb{P}_1[B_i].
\end{align*}
To bound the second term, we write 
\begin{align*}
\mathbb{P}_1[a_i = \theta,B_i] &= \mathbb{P}_1[a_i = \theta,a_i=a^*,B_i] + 
\mathbb{P}_1[a_i = \theta,a_i \neq a^*,B_i] \\
&= \mathbb{P}_1[a^* = \theta,a_i=a^*,B_i] + \mathbb{P}_1[a_i = \theta,a_i \neq
a^*,B_i]
\end{align*}
Since $\bar a$ satisfies herding in probability, $\lim_i\mathbb{P}%
[a_i=a^*]=1 $, in the limit the second term vanishes and it follows that 
\begin{align*}
\limsup_i\mathbb{P}_1[a_i = \theta,B_i] = \limsup_i \mathbb{P}_1[a^* =
\theta,B_i].
\end{align*}

It follows from the Concentration of Dependence Principle that
\begin{align*}
\limsup_i\mathbb{P}_1[a^* = \theta,B_i] = \mathbb{P}_1[a^* = \theta] \cdot \mathbb{P}_1[B_i] = p_1 \cdot \mathbb{P}_1[B_i],
\end{align*}
where the right-hand side does not depend on $i$, since the events $B_i$ all
have the same probability. We have thus shown that 
\begin{align*}
\limsup_i\mathbb{P}_1[a_i = \theta,B_i] = p_1 \cdot \mathbb{P}_1[B_i],
\end{align*}

Combining the bounds on the two terms we get that the expected gain in
utility conditioned on $\theta=1$ is 
\begin{align*}
\liminf_i \mathbb{P}_1[b_i=\theta, B_i] - \mathbb{P}_1[a_i = \theta,B_i] &\geq
(1-\beta-\varepsilon-p_1)\mathbb{P}_1[B_i].
\end{align*}
Since we assumed that $p_1 \leq 1-\beta-2\varepsilon$ we have that this is at
least $\varepsilon\mathbb{P}_1[B_i]$, and in particular positive. This completes the proof of Proposition~\ref{prop:learn-to-agree-bounded}. 

\section{Proof of Theorem~\ref{thm:herding}}

In this section we prove Theorem~\ref{thm:herding}.

Since the SLE is weakly ordered, we can identify the agents with the set of natural numbers $\{1,2,\ldots\}$ in such a way that if $i > j$ then $i$ knows $j$'s action. Let 
\begin{align*}
x_i = \mathbb{P}[\theta = 1 \mid a_1,\ldots,a_i]
\end{align*}
be the sequence of \emph{public beliefs}, and let $q_i = \mathbb{P}[\theta = 1 \mid \ell_i, s_i]$
be agent $i$'s equilibrium belief. Note that, since each agent $i$ knows $%
\{a_1,\ldots,a_i\}$, 
\begin{align}  \label{eq:p_i-q_i}
x_i = \mathbb{E}\left[q_i \mid a_1,\ldots,a_i\right],
\end{align}
by the law of total expectations.

Note also that the action $1$ is optimal for beliefs $1/2$ and higher, and the
action $0$ is optimal for beliefs $1/2$ and lower. Therefore, and since $%
\bar a$ is an equilibrium, 
\begin{align}  \label{eq:f-q}
a_i=1 \Rightarrow q_i\geq 1/2\quad\text{ and }\quad a_i=0 \Rightarrow q_i\leq 1/2
\end{align}
and 
\begin{align}  \label{eq:f-q-max}
\mathbb{E}[u(a_i,\theta)\mid q_i] = \mathbb{P}[a_i=\theta \mid q_i] =
\max\{q_i,1-q_i\}.
\end{align}

We start with two simple claims regarding $a_i$ and $x_i$.

\begin{claim}
\label{clm:overturning} If $a_i = 1$ then $x_i \geq 1/2$. If $a_i = 0$ then $%
x_i \leq 1/2$.
\end{claim}

\begin{proof}
By~\eqref{eq:f-q} we have that $q_i \geq 1/2$ conditioned on  $a_i=1$.
Hence, by~\eqref{eq:p_i-q_i}, $x_i \geq 1/2$ conditioned on  $a_i=1$. An
analogous argument holds for the case $a_i=0$.
\end{proof}

\begin{claim}
\label{clm:maxp} $\mathbb{P}[a_i=\theta\mid x_i] = \max\{x_i,1-x_i\}$.
\end{claim}

\begin{proof}
By Claim~\ref{clm:overturning}  
\begin{align*}
\mathbb{P}[\theta=a_i\mid x_i] &= 
\begin{cases}
\mathbb{P}[\theta=1\mid x_i] & \text{if }x_i > 1/2 \\ 
\mathbb{P}[\theta=0\mid x_i] & \text{if }x_i < 1/2 \\ 
\mathbb{P}[\theta=a_i\mid x_i] & \text{if }x_i = 1/2.%
\end{cases}%
\end{align*}
By~\eqref{eq:p_i-q_i} and~\eqref{eq:f-q}, if $x_i=1/2$ then  $x_i=q_i$.
Therefore, and since  $\mathbb{P}[\theta=1\mid x_i] = x_i$, and  $\mathbb{P}[a_i=\theta\mid q_i=1/2]=1/2$ by~\eqref{eq:f-q-max},  
\begin{align*}
\mathbb{P}[\theta=a_i\mid x_i] &= 
\begin{cases}
x_i & \text{if }x_i > 1/2 \\ 
1-x_i & \text{if }x_i < 1/2 \\ 
1/2 & \text{if }x_i = 1/2.%
\end{cases}%
\end{align*}
Thus $\mathbb{P}[\theta=a_i\mid x_i] = \max\{x_i,1-x_i\}$.
\end{proof}

Let $x = \mathbb{P}[\theta = 1 \mid \bar a]$,
and note that $x_i$ is a bounded martingale that converges a.s.\ to $x$. It
thus follows from Claim~\ref{clm:overturning} that conditioned on $a_i$ taking both values infinitely often it holds that $x=1/2$. Thus, to prove
our theorem, we will show that the probability of $x=1/2$ is zero.
Accordingly, define the event 
\begin{align*}
F^0 = \{x = 1/2\},
\end{align*}
and for $\varepsilon > 0$ define the events 
\begin{align*}
F_i^{\varepsilon} = \{x_i \in (1/2-\varepsilon,1/2+\varepsilon)\}.
\end{align*}
The event $F^{\varepsilon}_i$ is the event that the public belief $x_i$ is
close to $1/2$. Since the sequence $(x_i)_i$ converges a.s.\ to $x$, we have
that 
\begin{align}  \label{eq:F-eps1}
\lim_{i \to \infty} \mathbb{P}[F_0 \setminus F_{\varepsilon}^i]=0
\end{align}
for every $\varepsilon > 0$, and that 
\begin{align}  \label{eq:F-eps2}
\lim_{\varepsilon\to 0}\limsup_{i \to \infty} \mathbb{P}[F^{\varepsilon}_i]\geq 
\mathbb{P}[F^0].
\end{align}
Thus, to prove that $\mathbb{P}[F^0]=0$---which, as we explained above,
proves the claim---it suffices to show that the left hand side of \eqref{eq:F-eps2} vanishes.

To this end, as in the proof of Theorem~\ref{thm:consensus}, let 
\begin{align*}
b_x = b(s_x) \in \operatornamewithlimits{argmax}_{a \in A}\mathbb{P}%
[\theta=a \mid F^0, s_x]
\end{align*}
be an optimal action  given an additional agent $x$'s private signal only, conditioned on $F^0$. Let $b_i = b(s_i)$. Note that
\begin{align*}
    \mathbb{P}[b_i = \theta] = \mathbb{P}[b_x = \theta] = \mathbb{P}[b_x = \theta \mid F^0] > \frac{1}{2},
\end{align*}
where the first equality follows from the fact that $b_i$ and $b_x$ are conditionally identically distributed, and the second from the fact that $b_x$ is conditionally independent of $F^0$. The inequality follows because private signals are informative.

Consider the deviation in which player $i$ chooses $b_i$ instead of $a_i$, whenever $F_i^{\varepsilon}$ occurs; this is possible, since $F_i^{\varepsilon}$ is $\sigma(\ell_i,s_i)$-measurable. Then player $i$'s gain in expected utility from
this deviation is 
\begin{align*}
\mathbb{P}[b_i = \theta, F^{\varepsilon}_i] - \mathbb{P}[a_i = \theta,
F^{\varepsilon}_i].
\end{align*}
We prove that the left-hand side of \eqref{eq:F-eps2} vanishes by showing
that if it does not then
\begin{align*}
\lim_{\varepsilon \to 0}\limsup_{i \to \infty} \mathbb{P}[b_i=\theta ,
F^{\varepsilon}_i] - \mathbb{P}[a_i = \theta, F^{\varepsilon}_i] > 0,
\end{align*}
and thus this is a profitable deviation for some $\varepsilon$ small enough
and $i$ large enough, contradicting the assumption that $\bar a$ is an SLE.

To this end, we note that 
\begin{align*}
\mathbb{P}[b_i = \theta, F^{\varepsilon}_i] \geq \mathbb{P}[b_i = \theta,
F^0] - \mathbb{P}[b_i = \theta, F^0 \setminus F^{\varepsilon}_i],
\end{align*}
since 
\begin{align*}
F^0 \setminus (F^0\setminus F^{\varepsilon}_i) = F^0 \cap
F^{\varepsilon}_i\subseteq F^{\varepsilon}_i.
\end{align*}
It thus follows by~\eqref{eq:F-eps1} that 
\begin{align*}
\liminf_{i \to \infty}\mathbb{P}[b_i = \theta, F^{\varepsilon}_i] \geq
\liminf_{i \to \infty}\mathbb{P}[b_i = \theta, F^0].
\end{align*}

By the Concentration of Dependence Principle, as used in Lemma~\ref{lem:apply-cod}, 
\begin{align*}
\lim_{i \to \infty}\mathbb{P}[b_i = \theta, F^0] = \mathbb{P}[b_i =
\theta]\cdot \mathbb{P}[F^0].
\end{align*}
Since private signals are informative, $\mathbb{P}[b_i =
\theta] > 1/2$, and so we have that 
\begin{align}  \label{eq:liminf}
\liminf_{\varepsilon \to 0}\liminf_{i \to \infty}\mathbb{P}[b_i = \theta,
F^{\varepsilon}_i] > \frac{1}{2}\mathbb{P}[F^0].
\end{align}

Now, 
\begin{align*}
\mathbb{P}[a_i=\theta \mid F^{\varepsilon}_i] &= \mathbb{E}\big[\mathbb{P}%
[a_i=\theta \mid x_i]\mid F^{\varepsilon}_i\big] = \mathbb{E}\big[\max\{x_i,1-x_i\} \mid F^{\varepsilon}_i\big],
\end{align*}
where the second equality is an application of Claim~\ref{clm:maxp}. Since $%
x_i \in (1/2-\varepsilon,1/2+\varepsilon)$ conditioned on $F^{\varepsilon}_i$%
, we get that 
\begin{align*}
\mathbb{P}[a_i=\theta, F^{\varepsilon}_i] < \left(\frac{1}{2}%
+\varepsilon\right)\cdot \mathbb{P}[F^{\varepsilon}_i].
\end{align*}
Therefore, by~\eqref{eq:F-eps2}, 
\begin{align*}
\lim_{\varepsilon \to 0}\limsup_{i \to \infty} \mathbb{P}[a_i=\theta ,
F^{\varepsilon}_i] \leq \frac{1}{2}\cdot\mathbb{P}[F^0].
\end{align*}

Therefore, in combination with ~\eqref{eq:liminf}, the expected profit from
deviating from $a_{i}$ to $b_{i}$ on $F_{i}^{\varepsilon }$ satisfies 
\begin{equation*}
\lim_{\varepsilon \rightarrow 0}\limsup_{i\rightarrow \infty }\mathbb{P}%
[b_{i}=\theta ,F_{i}^{\varepsilon }]-\mathbb{P}[a_{i}=\theta
,F_{i}^{\varepsilon }]>0,
\end{equation*}%
and thus this is a profitable deviation for some $\varepsilon $ small enough
and $i$ large enough. Hence it follows that $F^{0}$ has probability zero,
concluding the proof of Theorem~\ref{thm:herding}.

\section{Proof of Theorem~\ref{thm:co-finite}}

Identify the set of agents with the natural numbers $\{1,2,\ldots\}$. Since there are only finite many actions (in fact, two), by compactness there is some subset $(i_k)_k$ of the agents whose actions converge in probability to some random action $a^*$:
\begin{align}
\label{eq:a-i-k-a}
    \lim_{k \to \infty}\mathbb{P}[a_{i_k}=a^*]=1.
\end{align}
Assume towards a contradiction that there is another subset $(j_k)_k$ that converges to a different action. That is, assume that there is some random $b^*$ such that $\mathbb{P}[a^* \neq b^*] >0$ and 
\begin{align}
\label{eq:a-j-k-a}
    \lim_{k \to \infty}\mathbb{P}[a_{j_k}=b^*]=1.
\end{align}
It follows from \eqref{eq:a-i-k-a} and \eqref{eq:a-j-k-a} that
\begin{align}
\label{eq:a-i-a}
    \lim_{k \to \infty}\mathbb{P}[a_{i_k}=\theta]=\mathbb{P}[a^*=\theta]~~\text{and}~~\lim_{k \to \infty}\mathbb{P}[a_{j_k}=\theta]=\mathbb{P}[b^*=\theta].
\end{align}

By the equilibrium property, if agent $i$ observes $j$'s action (i.e., if $a_j$ is $\sigma(\ell_i,s_i)$-measurable), then $\mathbb{P}[a_i=\theta] \geq \mathbb{P}[a_j=\theta]$.
Since the SLE is almost weakly ordered, we have that for a fixed $j$ this indeed holds for all $i$ large enough, and so, taking the limit along the sequence $(i_k)_k$ yields by \eqref{eq:a-i-a} that $\mathbb{P}[a^*=\theta] \geq \mathbb{P}[a_j=\theta]$.
Taking now the limit along $(j_k)_k$ and applying \eqref{eq:a-i-a} again yields
\begin{align*}
    \mathbb{P}[a^*=\theta] \geq \mathbb{P}[b^*=\theta].
\end{align*}
By symmetry, we have that, in fact, this holds with equality.

Now,
\begin{align*}
\mathbb{P}[a^*=\theta] 
&= \mathbb{P}[a^*=\theta,a^*=b^*] +\mathbb{P}[a^*=\theta,a^* \neq b^*]     \\
&= \mathbb{P}[a^*=\theta,b^*=\theta,a^*=b^*] +\mathbb{P}[a^*=\theta,a^* \neq b^*],
\end{align*}
and likewise
\begin{align*}
\mathbb{P}[b^*=\theta] 
&= \mathbb{P}[a^*=\theta,b^*=\theta,a^*=b^*] +\mathbb{P}[b^*=\theta,a^* \neq b^*].
\end{align*}
Since $\mathbb{P}[b^*=\theta] = \mathbb{P}[a^*=\theta]$, subtracting these equations yields
\begin{align*}
    \mathbb{P}[a^*=\theta,a^* \neq b^*] = \mathbb{P}[b^*=\theta,a^* \neq b^*].
\end{align*}
Hence it must be that
\begin{align}
    \label{eq:a*b*}
    \mathbb{P}[b^*=\theta \mid a^* \neq b^*] = \mathbb{P}[a^*=\theta \mid a^* \neq b^*] = \frac{1}{2}.
\end{align}
Hence by Bayes' Law 
\begin{align}
    \label{eq:a*b*theta}
    \mathbb{P}[\theta=1\mid a^* \neq b^*] = \frac{1}{2}.
\end{align}

Let $F$ be the event $a^* \neq b^*$. As in the proof of Theorem~\ref{thm:consensus}, choose 
\begin{align*}
  b_x  =b(s_x)\in \operatornamewithlimits{argmax}_{a \in A}\mathbb{P}[\theta=a \mid s_i, F],
\end{align*}
and $b_i = b(s_i)$. As in the proof of Theorem~\ref{thm:herding}, we note that $\mathbb{P}[b_i=\theta]>1/2$, since signals are informative.

Let $F_i$ be the event that conditioned on $i$'s information, the probability that $a^* \neq b^*$ is at least $1/2$:
\begin{align*}
    F_{i} = \{\mathbb{P}[a^* \neq b^*\mid \ell_i,s_i] \geq 1/2\}.
\end{align*}
Since the SLE is almost weakly ordered, it follows that 
\begin{align}
    \label{eq:F-i-F}
    \lim_{i \to \infty}\mathbb{P}[F_i \triangle F] = 0,
\end{align}
That is, the probability that $F$ occurs but $F_i$ does not---or vice versa---is very small for large $i$; for large $i$, agent $i$ approximately knows if $a^* \neq b^*$. This is implied by almost weak ordering (and the Martingale Convergence Theorem), which implies that for each $j$ it holds for all $i$ large enough that agent $i$ knows $(a_1,\ldots,a_j)$.

Consider the strategy in which agent $i$ chooses $b_i$ whenever $F_i$ occurs, and otherwise plays $a_i$. Then her expected gain from this deviation is
\begin{align*}
    P_i = \mathbb{P}[b_i=\theta,F_i] -  \mathbb{P}[a_i=\theta,F_{i}]
\end{align*}
since on the event $F_{i}^c$ her utility is the same as when she does not deviate (as on this event she indeed plays $a_i$ and does not deviate).

Recalling that $F$ is the event that $a^* \neq b^*$, it follows from \eqref{eq:a*b*} and \eqref{eq:F-i-F} that
\begin{align*}
    \limsup_{i \to \infty} P_i 
    &= \limsup_{i \to \infty} \mathbb{P}[b_i=\theta,F_{i}] -  \mathbb{P}[a_i=\theta,F_{i}]\\
    &= \limsup_{i \to \infty} \mathbb{P}[b_i=\theta,F] -  \mathbb{P}[a_i=\theta,F]\\
    &= \limsup_{i \to \infty} \mathbb{P}[b_i=\theta,F] - \frac{1}{2}\mathbb{P}[F].
\end{align*}
By the Concentration of Dependence Principle, as applied in Lemma~\ref{lem:apply-cod},
\begin{align*}
    \lim_{i \to \infty}\mathbb{P}[b_i=\theta,F] = \mathbb{P}[b_i=\theta]\cdot\mathbb{P}[F].
\end{align*}
Hence
\begin{align*}
    \lim_{i \to \infty}P_i = \mathbb{P}[F]\cdot\left(\mathbb{P}[b_i =\theta]-\frac{1}{2}\right).
\end{align*}
Since signals are informative this is positive, and so we have reached a contradiction. This completes the proof of Theorem~\ref{thm:co-finite}.

\section{Proof of Theorem~\protect\ref{thm:game-equilibria}}

This proof is essentially a recasting of the proof of Proposition 2.1 in 
\cite{rsv} to our language and notation.

Fix an agent $i$. The case that $\delta=0$ or $T_i$ is finite is immediate.
We thus assume henceforth that $\delta>0$ and $|T_i|=\infty$.

Let 
\begin{align*}
v_i = \max_{a \in A}\mathbb{E}[u(a,\theta) \mid \mathfrak{k}_i,s_i]
\end{align*}
be the maximum expected utility agent $i$ can guarantee given what she
(asymptotically) knows at the end of the game.

Fix $(\mathfrak{k}_i,s_i)$ and $\varepsilon >0$, and let $\overline{A}_\varepsilon,\underline{A}_\varepsilon \subseteq A$
be the sets of actions given by 
\begin{align*}
\overline{A}_\varepsilon = \big\{a \in A \,:\ \mathbb{E}[u(a,\theta) \mid \mathfrak{k}_i,s_i] >
v_i-\varepsilon\big\}
\end{align*}
and 
\begin{align*}
\underline{A}_\varepsilon = \big\{b \in A \,:\ \mathbb{E}[u(b,\theta) \mid \mathfrak{k}_i,s_i] <
v_i-3\varepsilon\big\}.
\end{align*}
That is, $\overline{A}_\varepsilon$ is the set of actions that is $\varepsilon$-optimal, and $\underline{A}_\varepsilon$ is
the set of actions that is $3\varepsilon$-suboptimal---conditioned on the
information available to the player at the end of the game.

For $t \in T_i$ let
\begin{align*}
    \overline{U}_{\varepsilon,t} = \inf_{a \in \overline{A}_\varepsilon}\mathbb{E}[u(a,\theta)\mid \mathfrak{k}_{i,t},s_i]
\end{align*}
be the worst expected utility (given the information available to $i$ at time $t$) of any $a \in \overline{A}_\varepsilon$. Note that $\overline{U_{\varepsilon,t}}$ is a bounded supermartingale, and likewise
\begin{align*}
    \underline{U}_{\varepsilon,t} = \sup_{b \in \underline{A}_\varepsilon}\mathbb{E}[u(b,\theta)\mid \mathfrak{k}_{i,t},s_i],
\end{align*}
is a bounded submartingale, and hence both converge. Furthermore,
\begin{align*}
    \lim_{t \in T_i} \overline{U}_{\varepsilon,t} \geq v_i -\varepsilon\quad\quad\text{and}\quad\quad\lim_{t \in T_i} \underline{U}_{\varepsilon,t} \leq v_i -3\varepsilon.
\end{align*}
That is, for large enough $t$, agent $i$ will assign high expected utility  to {\em all} actions in $\overline{A}_\varepsilon$, and low expected utility to {\em all} actions in $\underline{A}_\varepsilon$. In particular, it will almost surely hold for all $t \in T_i$ large enough that
\begin{align*}
    \inf_{a \in \overline{A}_\varepsilon}\mathbb{E}[u(a,\theta)\mid \mathfrak{k}_{i,t},s_i] > \sup_{b \in \underline{A}_\varepsilon}\mathbb{E}[u(b,\theta)\mid \mathfrak{k}_{i,t},s_i]+\varepsilon,
\end{align*}
so that for all $t$ large enough agent $i$ will 
have a larger expected utility for any $a \in \overline{A}_\varepsilon$, as compared to any $b \in \underline{A}_\varepsilon$. Hence agent $i$ will eventually only choose actions in $\underline{A}_\varepsilon^c$, the complement of $\underline{A}_\varepsilon$. Therefore, and since $\underline{A}_\varepsilon^c$ is compact, any limit point of the sequence of actions of $i$ must be in $\underline{A}_\varepsilon^c$. Since this holds for all $\varepsilon>0$,  we have shown that every limit point of the actions of $i$ must be in $\cap_{\varepsilon>0}\underline{A}_\varepsilon^c$, which is equal to the set of actions that yields an expected utility $v_i$, conditioned on  $(\mathfrak{k}_i,s_i)$. This concludes the proof.

\section{Proof of Proposition~\ref{prop:mst-agreement}}
Let $D$ be the event that $\bar A_i=\{0,1\}$ for all $i$. That is, $D$ is the event that every agent chooses both actions infinitely often. By \cite[Theorem 5.1]{mst}, $D$ is equivalent to the event that $\bar A_i=\{0,1\}$ for {\em some} agent $i$; that is, if $\bar A_i=\{0,1\}$ for some $i$ then the same holds for all. Hence the event $D$ is $\sigma(\mathfrak{k}_i,s_i)$-measurable for every player $i$, since at the end of the game each player knows if she took both actions infinitely often or not.

Let $(\bar\ell,\bar a)$ be given by $\ell_i = \mathfrak{k}_i$ and let $a_i$ be equal to some (measurable) choice from $\bar A_i$. Then $(\bar\ell_i,\bar a)$ is an SLE, by Theorem~\ref{thm:game-equilibria}.  Also, $D$ is $\sigma(\ell_i,s_i)$-measurable. Since conditioned on $D$ the probability that $\theta=1$ is $1/2$, we can proceed as in the proof of Theorem~\ref{thm:consensus}. Assume by contradiction that $D$ has positive probability, and consider the deviation in which, whenever $D$ occurs, agent $i$, instead of choosing $a_i$, chooses $b_i=b(s_i)$, where $b$ satisfies
\begin{align*}
  b(s_x)\in \operatornamewithlimits{argmax}_{a \in A}\mathbb{P}[\theta=a \mid D, s_x].
\end{align*}
The same argument of the proof of Theorem~\ref{thm:consensus} shows that this is a profitable deviation, thus showing that the assumption that $D$ has positive probability leads to a contradiction.


\newpage
\begin{center}
    {\huge Social Learning Equilibria \\ Supplementary Material}
\end{center}
\section{Additional Extensions}
This supplementary appendix presents additional extensions and results. The first concerns the case of heterogeneous types of agents, with the corresponding result following immediately from the results established in the paper. The second extension relaxes the assumption of binary states and actions inherent in the canonical setting. 

\subsection{Heterogeneous preferences and priors}

We relax the homogeneity assumption and consider agents who have
different utility functions and/or prior beliefs (with full support). Assume that all agents share the same belief regarding the conditional signal distributions but there are finitely many different types. Agent $i$'s type is determined by her full-support prior belief on the binary state space and her utility function for a binary action. We assume that for each utility function it is strictly preferable to match the action with the state than to mismatch.  

Assume that 
the agents' types are common knowledge. The following can be established by
making slight adjustments to our proofs.

\begin{proposition}
\label{thm:csle-heterogeneity} In a canonical setting with finitely many
commonly known types and where signals are unbounded, every CSLE satisfies
information aggregation.
\end{proposition}

If signals are unbounded then in a CSLE all agents agree on the same action,
and additionally this action is optimal. Thus unbounded signals overcome
heterogeneity in priors and payoffs. The result follows from the fact that types are commonly known and there exists at least one type with infinitely many agents. Proposition \ref{thm:csle-heterogeneity} is
interesting to view in light of Aumann's \textquotedblleft agreeing to
disagree\textquotedblright\ result \citeyearpar{aumann}. He showed that if
agents share a common prior then common knowledge of posteriors implies
agreement. Proposition \ref{thm:csle-heterogeneity} shows that if signals are
unbounded then (common) knowledge of actions implies that agreement and information aggregation hold
among an infinite group of agents, even if priors and utility functions
differ.
\subsection{Beyond the Canonical Setting: Many States and Actions}

Another extension of our results beyond the canonical setting is to settings with more than two states and more than two actions. In this section we consider social learning settings in which signals are still conditionally independent---as in the canonical setting---but the set of states can be of any (finite) size, as can the set of actions. We show that our agreement result for CSLEs still holds, under an additional condition on the structure of the utility function and the private signals; this condition rules out some pathological cases in which disagreement can arise.\footnote{Clearly, there can be disagreement in equilibrium when signals are completely uninformative. Likewise, when two actions yield the same utility in some state then agents who learn the state can choose different actions in equilibrium. These issues are avoided in the canonical settings, where each state has a different uniquely optimal action and where signals are informative. In this more general setting a more complicated assumption is required.} 

We consider  a social learning setting $(N,A,\Theta ,u,S,\mu)$ with $N$ countably infinite, $A$ and $\Theta$ finite, and conditionally independent private signals; we will refer to this as a {\em finite} setting. We say that private signals are {\em always useful} if, for any prior $p \in \Delta(\Theta)$ for which more than one action maximizes expected utility, observing a conditionally independent private signal (distributed as the agents' signals $s_i$) strictly increases the expected utility of a rational agent.\footnote{\cite{amf} show that in a general sequential social learning game, in every equilibrium signals are never  useful at the limit belief.} That is, if $s_x$ is an additional conditionally independent private signal, distributed as $s_i$, and if $D$ is any event that is conditionally independent of $s_x$, then whenever 
$$
  \left|\argmax_{a \in A}\mathbb{E}[u(a,\theta)\mid D]\right| \geq 2 
$$
(that is, whenever conditioning on $D$ results in more than one optimal action) it holds that
$$
  \mathbb{E}[\max_{a \in A}\mathbb{E}[u(a,\theta)\mid D,s_x] \mid D] > \max_{a \in A}\mathbb{E}[u(a,\theta)\mid D]
$$
(that is, learning $s_x$ increases one's  expected utility).
It is easy to see that for the case of two states and two actions, this holds whenever signals are informative.

The assumption of always useful signals implies that it is impossible for a state $\theta_0$ to have more than one optimal action, since otherwise, conditioning on the event that the state is $\theta_0$, an additional private signal will not change the agent's belief and thus cannot result in higher expected utility. Another implication for the case of many states is that signals cannot be restricted to some particular dimension and ignore others. For example, if $\Theta = \{0,1\}\times\{0,1\}$, the assumption of always useful signals rules out a signal that is informative with respect to the first coordinate, but provides no information regarding the second. Note that the assumption of always useful signals does not imply that signals are unbounded.

\begin{proposition}
\label{thm:clse_agreement_many_states_actions}
Consider a finite setting with infinitely many agents. If signals are always useful then every CSLE satisfies agreement.
\end{proposition}
The information aggregation results of Proposition~\ref{prp:csle-aggregating} and Theorem~\ref{thm:learn-to-agree} also hold in this setting, under an appropriate definition of unbounded signals. In interest of brevity we leave Proposition~\ref{thm:clse_agreement_many_states_actions} as the only result that we extend in this direction.

The proof of Proposition~\ref{thm:clse_agreement_many_states_actions} starts with the following ``No Trade'' Lemma, showing that disagreement implies indifference. Its proof, which we omit, follows the same argument as the proof of Lemma~\ref{lem:no-bet}.
\begin{lemma}
\label{lem:no-bet-many} 
Fix a social learning setting  with $A$ and $\Theta$ finite, and let $(\bar \ell,\bar a)$ be an SLE defined on this setting. Let $A_0$ be a subset of $A$, and let $D$ be the event that for each $a \in A_0$ there is an agent who chooses the action $a$: 
\begin{align*}
  D = \{\text{for each}~a \in A_0~\text{there exists an}~i \in N~\text{such that}~a_i=a\}.
\end{align*}
If the probability of $D$ is positive, then for any $a,b \in A_0$ it holds that
\begin{align*}
    \mathbb{E}[u(a,\theta)\mid D] = \mathbb{E}[u(b,\theta)\mid D].
\end{align*}
\end{lemma}

Given this lemma, we turn to the proof of Proposition~\ref{thm:clse_agreement_many_states_actions}.
\begin{proof}[Proof of Proposition~\ref{thm:clse_agreement_many_states_actions}]
  Let $A_0$ be a subset of $A$ of size at least $2$, and let $D$ be the disagreement event that for each $a \in A_0$ there is an agent who chooses the action $a$: 
  \begin{align*}
      D = \{\text{for each}~a \in A_0~\text{there exists an}~i \in N~\text{such that}~a_i=a\}.
  \end{align*}
  We assume by contradiction that $D$ has positive probability. By Lemma~\ref{lem:no-bet-many},
  \begin{align*}
      \mathbb{E}[u(a,\theta)\mid D] = \mathbb{E}[u(b,\theta)\mid D]
  \end{align*}
  for all $a,b \in A_0$, and we denote this quantity by $U(A_0\mid D)$. Using this notation, we can write the expected utility of agent $i$  as
  \begin{align}
  \label{eq:thm_1_no_dev}
      \mathbb{E}[u(a_i,\theta)\ind{D}] + \mathbb{E}[u(a_i,\theta)(1-\ind{D})]
      &= U(A_0 \mid D) \cdot \mathbb{P}[D] + \mathbb{E}[u(a_i,\theta)(1-\ind{D})]
  \end{align}

  As in the proof of Theorem~\ref{thm:consensus}, let
  \begin{align*}
    b_x  =b(s_x) \in \operatornamewithlimits{argmax}_{a \in A_0}\mathbb{E}[u(a,\theta)\mid D,s_x],
  \end{align*}
  denote $b_i = b(s_i)$, and consider a deviation by agent $i$ in which she chooses $b_i$ whenever $D$ occurs, and $a_i$ otherwise. Then the expected utility of this deviation is
  \begin{align}
     \label{eq:thm_1_dev}
      \mathbb{E}[u(b_i,\theta)\ind{D}] + \mathbb{E}[u(a_i,\theta)(1-\ind{D})].
  \end{align}
  Now,
  \begin{align}
  \label{eq:lim_exp_utl}
      \mathbb{E}[u(b_i,\theta)\ind{D}] 
      &= \sum_{\omega \in \Theta} \mathbb{E}[u(b_i,\theta)\ind{D}\mid \theta=\omega]\cdot\mathbb{P}[\theta=\omega]\nonumber\\
      &= \sum_{\omega \in \Theta} \mathbb{E}[u(b_i,\omega)\ind{D}\mid \theta=\omega]\cdot\mathbb{P}[\theta=\omega]\nonumber\\
      &= \sum_{\omega \in \Theta} \sum_{a \in A}u(a,\omega)\mathbb{P}[b_i=a,D\mid \theta=\omega]\cdot\mathbb{P}[\theta=\omega].
  \end{align}
  By the Concentration of Dependence Principle 
  \begin{align*}
      \lim_{i \to \infty}\mathbb{P}[b_i=a,D\mid \theta=\omega] = \mathbb{P}[b_i=a\mid \theta=\omega] \cdot \mathbb{P}[D\mid \theta=\omega].
  \end{align*}
  Note that the right hand side holds for any $i$, since signals are conditionally identically distributed. Hence it also holds for agent $x$, and thus
  \begin{align*}
      \lim_{i \to \infty}\mathbb{P}[b_i=a,D\mid \theta=\omega] = \mathbb{P}[b_x=a\mid \theta=\omega] \cdot \mathbb{P}[D\mid \theta=\omega].
  \end{align*}
  Substituting this back into \eqref{eq:lim_exp_utl} yields
  \begin{align*}
      \lim_{i \to \infty}\mathbb{E}[u(b_i,\theta)\ind{D}] 
      &= \sum_{\omega \in \Theta} \sum_{a \in A}u(a,\omega)\mathbb{P}[b_x=a\mid \theta=\omega]\cdot\mathbb{P}[D\mid \theta=\omega]\cdot\mathbb{P}[\theta=\omega]\\
      &= \sum_{\omega \in \Theta} \sum_{a \in A}u(a,\omega)\mathbb{P}[b_x=a,D\mid \theta=\omega]\cdot\mathbb{P}[\theta=\omega]\\
      &= \mathbb{E}[u(b_x,\theta)\mid D] \cdot \mathbb{P}[D].
  \end{align*}
  By our assumption that the signals are always useful, 
  \begin{align*}
      \mathbb{E}[u(b_x,\theta)\mid D] > U(A_0\mid D).
  \end{align*}
  Substituting this back into \eqref{eq:thm_1_dev} and comparing to \eqref{eq:thm_1_no_dev} shows that $b_i$ is a profitable deviation for some $i$ large enough, and so we have reached a contradiction with our equilibrium assumption.
\end{proof}

\end{document}